\documentclass[12pt]{article}

\newcommand{\figPremierePage}{
	\begin{tikzpicture}[node distance=1.8cm,on grid,scale=1.5]
		\pgfmathtruncatemacro{\n}{8}
		\pgfmathtruncatemacro{\nUn}{\n-1}
		\pgfmathtruncatemacro{\nDeux}{\n-2}
		\pgfmathtruncatemacro{\nTrois}{\n-3}
		\begin{scope}[name prefix = b-]
			\node (1) at (0,0) {};
			\foreach \x in {2,...,\n}
			{
				\pgfmathtruncatemacro{\prec}{\x - 1}
				\node (\x) [right=of \prec] {};
			}
		\end{scope}
		\begin{scope}[name prefix = t-]
			\node (1) at ($(b-1)+(0,4)$) {};
			\foreach \x in {2,...,\n}
			{
				\pgfmathtruncatemacro{\prec}{\x - 1}
				\node (\x) [right=of \prec] {};
			}
		\end{scope}
		
		\path[{Circle[]}-{Circle[]}]
		(t-1) edge [bend right=90,looseness=2] (t-2)
		(t-5) edge [bend right=90,looseness=2] (t-6)
		(t-4) edge [bend right=90,looseness=1.6] (t-7)
		
		(b-1) edge [bend left=90,looseness=2] (b-2)
		(b-4) edge [bend left=90,looseness=2] (b-5)
		(b-7) edge [bend left=90,looseness=2] (b-8)
		
		(b-6) edge [out=90,in=-90] (t-8)
		
		(b-3) edge  (t-3);
	\end{tikzpicture}
}

\newcommand{\figQuatreDiagramme}{
	\begin{tikzpicture}[node distance=0.5cm]
		\begin{scope}[name prefix = t]
			\node (1) at (0,2) {1};
			\node (2) [right=of 1] {2};
			\node (3) [right=of 2] {3};
			\node (4) [right=of 3] {4};
		\end{scope}
		\begin{scope}[name prefix = b]
			\node (1) at (0,0) {1};
			\node (2) [right=of 1] {2};
			\node (3) [right=of 2] {3};
			\node (4) [right=of 3] {4};
		\end{scope}
		\path[{Circle[]}-{Circle[]}]   (t1) edge (b1)
		(t2) edge [bend right=90,looseness=2] (t3)
		(t4) edge [out=-90,in=90] (b2)
		(b3) edge [bend left=90,looseness=2] (b4);
	\end{tikzpicture}
}

\newcommand{\figDelta}{
	\begin{tikzpicture}[node distance=1cm,on grid]
		\begin{scope}[name prefix = ml]
			\node (1) at (0,0) {};
			\node (2) [right=of 1] {};
			\node (3) [right=of 2] {};
		\end{scope}
		\begin{scope}[name prefix = tl]
			\node (1) at ($(ml1)+(0,2)$) {1};
			\node (2) [right=of 1] {2};
			\node (3) [right=of 2] {3};
		\end{scope}
		\begin{scope}[name prefix = bl]
			\node (1) at ($(ml1)+(0,-2)$) {1};
			\node (2) [right=of 1] {2};
			\node (3) [right=of 2] {3};
		\end{scope}
		
		\node (=) at ($(ml3)+(1,0)$) {\huge{$= \delta$}};
		
		\begin{scope}[name prefix = mr]
			\node (1) at ($(ml3)+(2,0)$) {};
			\node (2) [right=of 1] {};
			\node (3) [right=of 2] {};
		\end{scope}	
		\begin{scope}[name prefix = tr]
			\node (1) at ($(mr1)+(0,2)$) {1};
			\node (2) [right=of 1] {2};
			\node (3) [right=of 2] {3};
		\end{scope}
		\begin{scope}[name prefix = br]
			\node (1) at ($(mr1)+(0,-2)$) {1};
			\node (2) [right=of 1] {2};
			\node (3) [right=of 2] {3};
		\end{scope}
		
		\path[{Circle[]}-{Circle[]}] 
		(tl1) edge (ml1.center)  
		(tr1) edge (mr1.center)
		(tl2) edge [bend right=90,looseness=2] (tl3)
		(bl2) edge [bend left=90,looseness=2] (bl3)
		(ml2.center) edge [bend left=90,looseness=2] (ml3.center)
		(tr2) edge [bend right=90,looseness=2] (tr3)
		(br2) edge [bend left=90,looseness=2] (br3);
		\path[{Circle[]}-]   
		(bl1) edge (ml1.center)
		(br1) edge (mr1.center);
		\path 
		(ml3.center) edge [bend left=90,looseness=2] (ml2.center);
		
	\end{tikzpicture}
}

\newcommand{\figDiagrammeToPont}{
	\begin{tikzpicture}[node distance=1cm,on grid]
		\begin{scope}[name prefix = ml]
			\node (1) at (0,0) {1};
			\node (2) [right=of 1] {2};
			\node (3) [right=of 2] {3};
			\node (4) [right=of 3] {4};
		\end{scope}
		\begin{scope}[name prefix = tl]
			\node (1) at ($(ml1)+(0,2)$) {8};
			\node (2) [right=of 1] {7};
			\node (3) [right=of 2] {6};
			\node (4) [right=of 3] {5};
		\end{scope}
		
		\begin{scope}[name prefix = r]
			\node (1) at ($(ml4)+(2,0)$) {1};
			\node (2) [right=of 1] {2};
			\node (3) [right=of 2] {3};
			\node (4) [right=of 3] {4};
			\node (5) [right=of 4] {5};
			\node (6) [right=of 5] {6};
			\node (7) [right=of 6] {7};
			\node (8) [right=of 7] {8};
		\end{scope}
		
		\path[{Circle[]}-{Circle[]}]   
		(tl1) edge (ml1)
		(tl2) edge [bend right=90,looseness=2] (tl3)
		(tl4) edge [out=-90,in=90] (ml2)
		(ml3) edge [bend left=90,looseness=2] (ml4)
		
		(r1) edge [bend left=90,looseness=1] (r8)
		(r2) edge [bend left=90,looseness=1.5] (r5)
		(r3) edge [bend left=90,looseness=2] (r4)
		(r6) edge [bend left=90,looseness=2] (r7);
		
		\node[draw=blue,inner sep=3pt,ellipse,fit=(tl1) (tl2) (tl3) (tl4)] {};
		\node[draw=blue,inner sep=2pt,ellipse,fit=(r5) (r6) (r7) (r8)] {};
		\draw[blue,->] (1.5,2.5) to [bend left=70,looseness=0.5] (9,2.5);
	\end{tikzpicture}
}

\newcommand{\figPont}{
	\begin{tikzpicture}[node distance=1cm,on grid]
		
		\node (1) at (0,0) {1};
		\node (2) [right=of 1] {2};
		\node (3) [right=of 2] {3};
		\node (4) [right=of 3] {4};
		\node (5) [right=of 4] {5};
		\node (6) [right=of 5] {k};
		\node (7) [right=of 6] {7};
		\node (8) [right=of 7] {8};
		
		\draw[{Circle[]}-{Circle[]}] (1) edge [bend left=90,looseness=1] (6);
		
		\node (22) at (2.north)  {\huge{\textbf{.}}};
		\node (33) at (3.north)  {\huge{\textbf{.}}};
		\node (44) at (4.north)  {\huge{\textbf{.}}};
		\node (55) at (5.north)  {\huge{\textbf{.}}};
		\node (77) at (7.north)  {\huge{\textbf{.}}};
		\node (88) at (8.north)  {\huge{\textbf{.}}};
		
	\end{tikzpicture}
}

\newcommand{\figGenerateurs}{
	\begin{tikzpicture}[node distance=0.5cm,on grid]
		\begin{scope}[name prefix = m1-]
			\node (1) at (-20,0) {1};
			\node (2) [right=of 1] {2};
			\node (3) [right=of 2] {3};
			\node (4) [right=of 3] {4};
		\end{scope}
		\begin{scope}[name prefix = t1-]
			\node (1) at ($(m1-1)+(0,2)$) {1};
			\node (2) [right=of 1] {2};
			\node (3) [right=of 2] {3};
			\node (4) [right=of 3] {4};
		\end{scope}
		
		\path[{Circle[]}-{Circle[]}]   
		(m1-1) edge (t1-1)
		(m1-2) edge (t1-2)
		(m1-3) edge (t1-3)
		(m1-4) edge (t1-4);
		
		\begin{scope}[name prefix = m2-]
			\node (1) at ($(m1-4)+(2,0)$) {1};
			\node (2) [right=of 1] {2};
			\node (3) [right=of 2] {3};
			\node (4) [right=of 3] {4};
		\end{scope}
		\begin{scope}[name prefix = t2-]
			\node (1) at ($(m2-1)+(0,2)$) {1};
			\node (2) [right=of 1] {2};
			\node (3) [right=of 2] {3};
			\node (4) [right=of 3] {4};
		\end{scope}
		
		\path[{Circle[]}-{Circle[]}]   
		(t2-1) edge [bend right=90,looseness=2] (t2-2)
		(m2-1) edge [bend left=90,looseness=2] (m2-2)
		(m2-3) edge (t2-3)
		(m2-4) edge (t2-4);
		
		\begin{scope}[name prefix = m3-]
			\node (1) at ($(m2-4)+(2,0)$) {1};
			\node (2) [right=of 1] {2};
			\node (3) [right=of 2] {3};
			\node (4) [right=of 3] {4};
		\end{scope}
		\begin{scope}[name prefix = t3-]
			\node (1) at ($(m3-1)+(0,2)$) {1};
			\node (2) [right=of 1] {2};
			\node (3) [right=of 2] {3};
			\node (4) [right=of 3] {4};
		\end{scope}
		
		\path[{Circle[]}-{Circle[]}]   
		(t3-2) edge [bend right=90,looseness=2] (t3-3)
		(m3-2) edge [bend left=90,looseness=2] (m3-3)
		(m3-1) edge (t3-1)
		(m3-4) edge (t3-4);
		
		\begin{scope}[name prefix = m4-]
			\node (1) at ($(m3-4)+(2,0)$) {1};
			\node (2) [right=of 1] {2};
			\node (3) [right=of 2] {3};
			\node (4) [right=of 3] {4};
		\end{scope}
		\begin{scope}[name prefix = t4-]
			\node (1) at ($(m4-1)+(0,2)$) {1};
			\node (2) [right=of 1] {2};
			\node (3) [right=of 2] {3};
			\node (4) [right=of 3] {4};
		\end{scope}
		
		\path[{Circle[]}-{Circle[]}]   
		(t4-3) edge [bend right=90,looseness=2] (t4-4)
		(m4-3) edge [bend left=90,looseness=2] (m4-4)
		(m4-1) edge (t4-1)
		(m4-2) edge (t4-2);
		
	\end{tikzpicture}
}

\newcommand{\figCommutativite}{
	\begin{tikzpicture}[node distance=1cm,on grid,scale=0.8]
		\begin{scope}[name prefix = m2-]
			\node (1) at (0,0) {};
			\node (2) [right=of 1] {};
			\node (3) [right=of 2] {};
			\node (4) [right=of 3] {};
		\end{scope}
		\begin{scope}[name prefix = t2-]
			\node (1) at ($(m2-1)+(0,2)$) {1};
			\node (2) [right=of 1] {2};
			\node (3) [right=of 2] {3};
			\node (4) [right=of 3] {4};
		\end{scope}
		
		\begin{scope}[name prefix = m4-]
			\node (1) at ($(m2-1)+(0,-2)$) {1};
			\node (2) [right=of 1] {2};
			\node (3) [right=of 2] {3};
			\node (4) [right=of 3] {4};
		\end{scope}
		\begin{scope}[name prefix = t4-]
			\node (1) at ($(m4-1)+(0,2)$) {};
			\node (2) [right=of 1] {};
			\node (3) [right=of 2] {};
			\node (4) [right=of 3] {};
		\end{scope}
		
		\begin{scope}[name prefix = m3-]
			\node (1) at ($(m4-4)+(2,0)$) {1};
			\node (2) [right=of 1] {2};
			\node (3) [right=of 2] {3};
			\node (4) [right=of 3] {4};
		\end{scope}
		\begin{scope}[name prefix = t3-]
			\node (1) at ($(m3-1)+(0,2)$) {};
			\node (2) [right=of 1] {};
			\node (3) [right=of 2] {};
			\node (4) [right=of 3] {};
		\end{scope}
		
		\begin{scope}[name prefix = m1-]
			\node (1) at ($(m2-4)+(2,0)$) {};
			\node (2) [right=of 1] {};
			\node (3) [right=of 2] {};
			\node (4) [right=of 3] {};
		\end{scope}
		\begin{scope}[name prefix = t1-]
			\node (1) at ($(m1-1)+(0,2)$) {1};
			\node (2) [right=of 1] {2};
			\node (3) [right=of 2] {3};
			\node (4) [right=of 3] {4};
		\end{scope}
		
		\path[{Circle[]}-{Circle[]}]   
		(t2-1) edge [bend right=90,looseness=1.5] (t2-2)
		(m4-3) edge [bend left=90,looseness=1.5] (m4-4)
		(m3-1) edge [bend left=90,looseness=1.5] (m3-2)
		(t1-3) edge [bend right=90,looseness=1.5] (t1-4);
		
		\path[{Circle[]}-] 
		(t2-3) edge (m2-3.center)
		(t2-4) edge (m2-4.center)
		(m4-1) edge (t4-1.center) 
		(m4-2) edge (t4-2.center)
		(m3-3) edge (t3-3.center)
		(m3-4) edge (t3-4.center)
		(t1-1) edge (m1-1.center)
		(t1-2) edge (m1-2.center);
		
		\path
		(m2-1.center) edge [bend left=90,looseness=1.5] (m2-2.center)
		(t4-3.center) edge [bend right=90,looseness=1.5] (t4-4.center)
		(t3-1.center) edge [bend right=90,looseness=1.5] (t3-2.center)
		(m1-3.center) edge [bend left=90,looseness=1.5] (m1-4.center);
		
		\node at ($(m2-4)+(1,0)$) {\huge{\textbf{=}}};
		
	\end{tikzpicture}
}

\newcommand{\figProche}{
	\begin{tikzpicture}[node distance=1cm,on grid,scale=0.8]
		\begin{scope}[name prefix = m2-]
			\node (1) at (0,0) {};
			\node (2) [right=of 1] {};
			\node (3) [right=of 2] {};
			\node (4) [right=of 3] {};
		\end{scope}
		\begin{scope}[name prefix = t2-]
			\node (1) at ($(m2-1)+(0,2)$) {1};
			\node (2) [right=of 1] {2};
			\node (3) [right=of 2] {3};
			\node (4) [right=of 3] {4};
		\end{scope}
		
		\begin{scope}[name prefix = m4-]
			\node (1) at ($(m2-1)+(0,-2)$) {};
			\node (2) [right=of 1] {};
			\node (3) [right=of 2] {};
			\node (4) [right=of 3] {};
		\end{scope}
		\begin{scope}[name prefix = t4-]
			\node (1) at ($(m4-1)+(0,2)$) {};
			\node (2) [right=of 1] {};
			\node (3) [right=of 2] {};
			\node (4) [right=of 3] {};
		\end{scope}
		
		\begin{scope}[name prefix = m3-]
			\node (1) at ($(m4-4)+(2,-0.5)$) {1};
			\node (2) [right=of 1] {2};
			\node (3) [right=of 2] {3};
			\node (4) [right=of 3] {4};
		\end{scope}
		\begin{scope}[name prefix = t3-]
			\node (1) at ($(m3-1)+(0,3)$) {1};
			\node (2) [right=of 1] {2};
			\node (3) [right=of 2] {3};
			\node (4) [right=of 3] {4};
		\end{scope}
		
		\begin{scope}[name prefix = m1-]
			\node (1) at ($(m4-1)+(0,-2)$) {1};
			\node (2) [right=of 1] {2};
			\node (3) [right=of 2] {3};
			\node (4) [right=of 3] {4};
		\end{scope}
		\begin{scope}[name prefix = t1-]
			\node (1) at ($(m1-1)+(0,2)$) {};
			\node (2) [right=of 1] {};
			\node (3) [right=of 2] {};
			\node (4) [right=of 3] {};
		\end{scope}
		
		\path[{Circle[]}-{Circle[]}]  
		(m1-1) edge [bend left=90,looseness=1.5] (m1-2)
		(t2-1) edge [bend right=90,looseness=1.5] (t2-2) 
		(t3-1) edge [bend right=90,looseness=2] (t3-2)
		(m3-1) edge [bend left=90,looseness=2] (m3-2)
		(m3-3) edge (t3-3)
		(m3-4) edge (t3-4);
		
		\path[{Circle[]}-] 
		(m1-3) edge (t1-3.center)
		(m1-4) edge (t1-4.center)  
		(t2-3) edge (m2-3.center)
		(t2-4) edge (m2-4.center);
		
		\path
		(t1-1.center) edge [bend right=90,looseness=1.5] (t1-2.center)
		(m2-1.center) edge [bend left=90,looseness=1.5] (m2-2.center)
		(t4-2.center) edge [bend right=90,looseness=1.5] (t4-3.center)
		(m4-2.center) edge [bend left=90,looseness=1.5] (m4-3.center)
		(m4-1.center) edge (t4-1.center)
		(m4-4.center) edge (t4-4.center);
		
		\node at ($(m2-4)+(1,-1)$) {\huge{\textbf{=}}};
		
	\end{tikzpicture}
}

\newcommand{\figDiagNonCommutatif}{
	\begin{tikzpicture}[node distance=1.5cm]
		\node (ccsd) at (0,0) {$CCSD(n)$};
		\node (mn) [right=of ccsd] {$M_n/{\sim}$};
		\node (dn) [above=of ccsd] {$D_n$};
		
		\path[->]
		(dn) edge (ccsd) 
		(ccsd) edge (mn)
		(dn) edge node[above] {$\varphi$} (mn);	
	\end{tikzpicture}
}

\newcommand{\figCCSDEi}{
	\begin{tikzpicture}[node distance=1cm,on grid]
		\begin{scope}[name prefix = m2-]
			\node (1) at (0,0) {1};
			\node (2) [right=of 1] {2};
			\node (3) [right=of 2] {3};
			\node (4) [right=of 3] {4};
		\end{scope}
		\begin{scope}[name prefix = t2-]
			\node (1) at ($(m2-1)+(0,2)$) {1};
			\node (2) [right=of 1] {2};
			\node (3) [right=of 2] {3};
			\node (4) [right=of 3] {4};
		\end{scope}
		
		\path[{Circle[]}-{Circle[]}]   
		(t2-1) edge [bend right=90,looseness=2] (t2-2)
		(m2-1) edge [bend left=90,looseness=2] (m2-2)
		(m2-3) edge (t2-3)
		(m2-4) edge (t2-4);

		\begin{scope}[name prefix = r-]
			\node (1) at ($(m2-4)+(2,0)$) {1};
			\node (2) [right=of 1] {2};
			\node (3) [right=of 2] {3};
			\node (4) [right=of 3] {4};
			\node (5) [right=of 4] {5};
			\node (6) [right=of 5] {6};
			\node (7) [right=of 6] {7};
			\node (8) [right=of 7] {8};
		\end{scope}

		\path[{Circle[]}-{Circle[]}]   
		(r-1) edge [bend left=90,looseness=2] (r-2)
		(r-7) edge [bend left=90,looseness=2] (r-8)
		(r-3) edge [bend left=90,looseness=2] (r-6)
		(r-4) edge [bend left=90,looseness=2] (r-5);
		
		\node at ($(m2-4)+(1,1)$) {\huge{\textbf{$\rightarrow$}}};
		
	\end{tikzpicture}
	\begin{tikzpicture}[scale=1]
		\tikzset{
			grid_style/.style = {
				step=1cm,
				gray,
				very thin,
				dashed
			},
			point_style/.style = {
				draw,
				circle,
				inner sep=1pt
			},
			path_style/.style = {
				blue,
				thick,
				->,
				looseness=1
			}
		}
		
		\draw[grid_style] (0,0) grid (4,4);
		
		\foreach \x in {0,...,4}
		\foreach \y in {0,...,4}
		{
			\node[point_style] at (\x,\y) {};
		}
		
		\draw[path_style] (0,0) --++ (1,0) --++ (0,1) --++ (2,0)--++ (0,2) --++ (1,0) --++ (0,1);
		\draw[red] (0,0) -- (4,4);
		
		\foreach \x in {0,...,4}
		{
			\node at (\x,0) [below] {\x};
			\node at (0,\x) [left] {\x};
		}
	\end{tikzpicture}
}

\newcommand{\figCCSDei}{
	\begin{tikzpicture}[scale=1]
		\tikzset{
			grid_style/.style = {
				step=1cm,
				gray,
				very thin,
				dashed
			},
			point_style/.style = {
				draw,
				circle,
				inner sep=1pt
			},
			path_style/.style = {
				blue,
				thick,
				->,
				looseness=1
			}
		}
		
		\draw[grid_style] (0,0) grid (4,4);
		
		\foreach \x in {0,...,4}
		\foreach \y in {0,...,4}
		{
			\node[point_style] at (\x,\y) {};
		}
		
		\draw[path_style] (0,0) --++ (1,0) --++ (0,1) --++ (3,0)--++ (0,3);
		\draw[red] (0,0) -- (4,4);
		
		\foreach \x in {0,...,4}
		{
			\node at (\x,0) [below] {\x};
			\node at (0,\x) [left] {\x};
		}
	\end{tikzpicture}
}

\newcommand{\figSegmentFNJ}{
	\begin{tikzpicture}[node distance=1cm,on grid]
		\pgfmathtruncatemacro{\n}{4}
		\pgfmathtruncatemacro{\nUn}{\n-1}
		\pgfmathtruncatemacro{\nDeux}{\n-2}
		\pgfmathtruncatemacro{\nTrois}{\n-3}
		\begin{scope}[name prefix = m-]
			\node (1) at (0,0) {1};
			\foreach \x in {2,...,\n}
			{
				\pgfmathtruncatemacro{\prec}{\x - 1}
				\node (\x) [right=of \prec] {\x};
			}
		\end{scope}
		\begin{scope}[name prefix = t-]
			\node (1) at ($(m-1)+(0,3)$) {1};
			\foreach \x in {2,...,\n}
			{
				\pgfmathtruncatemacro{\prec}{\x - 1}
				\node (\x) [right=of \prec] {\x};
			}
		\end{scope}
		
		\path[{Circle[]}-{Circle[]}]
		(t-1) edge [bend right=90,looseness=2] (t-2)
		(m-\nUn) edge [bend left=90,looseness=2] (m-\n);
		
		\foreach \x in {1,...,\nDeux}
		{
			\pgfmathtruncatemacro{\suiv}{\x + 2}
			\draw[{Circle[]}-{Circle[]}] (m-\x) edge [out=90,in=-90] (t-\suiv);
		}
		
		\node at ($(m-1)+(-1,1.5)$) {\textbf{\huge{=}}};
		
		\begin{scope}[name prefix = t0-]
			\node (1) at ($(t-1)+(-5,1.5)$) {1};
			\foreach \x in {2,...,\n}
			{
				\pgfmathtruncatemacro{\prec}{\x - 1}
				\node (\x) [right=of \prec] {\x};
			}
		\end{scope}
		\foreach \i in {1,2}
		{
			\begin{scope}[name prefix =t\i-]
				\pgfmathtruncatemacro{\top}{\i - 1}
				\node (1) at ($(t\top-1)+(0,-2)$) {};
				\foreach \x in {2,...,\n}
				{
					\pgfmathtruncatemacro{\prec}{\x - 1}
					\node (\x) [right=of \prec] {};
				}
			\end{scope}
		}
		\begin{scope}[name prefix =t3-]
			\node (1) at ($(t2-1)+(0,-2)$) {1};
			\foreach \x in {2,...,\n}
			{
				\pgfmathtruncatemacro{\prec}{\x - 1}
				\node (\x) [right=of \prec] {\x};
			}
		\end{scope}
		
		\path[{Circle[]}-{Circle[]}] 
		(t0-1) edge[bend right = 90, looseness=2] (t0-2)
		(t3-3) edge[bend left = 90, looseness=2] (t3-4);
		
		\path[{Circle[]}-] 
		(t0-3) edge (t1-3.center)
		(t0-4) edge (t1-4.center)
		(t3-1) edge (t2-1.center)
		(t3-2) edge (t2-2.center);
		
		\path 
		(t1-1.center) edge (t2-1.center)
		(t1-4.center) edge (t2-4.center)
		(t1-1.center) edge[bend left = 90, looseness=2] (t1-2.center)
		(t1-2.center) edge[bend right = 90, looseness=2] (t1-3.center)
		(t2-2.center) edge[bend left = 90, looseness=2] (t2-3.center)
		(t2-3.center) edge[bend right = 90, looseness=2] (t2-4.center);

	\end{tikzpicture}
}

\newcommand{\figSeparationSegmentFNJ}{
	\begin{tikzpicture}[node distance=1cm,on grid]
		\pgfmathtruncatemacro{\n}{7}
		\pgfmathtruncatemacro{\nUn}{\n-1}
		\pgfmathtruncatemacro{\nDeux}{\n-2}
		\pgfmathtruncatemacro{\nTrois}{\n-3}
		
		\begin{scope}[name prefix = m-]
			\node (1) at (0,0) {};
			\foreach \x in {2,...,\n}
			{
				\pgfmathtruncatemacro{\prec}{\x - 1}
				\node (\x) [right=of \prec] {};
			}
			\node (l) [left=of 1] {};
			\node (0) [left=of l] {};
			\node (rdots)  [right=of \n] {};
			\node (n)  [right=of rdots] {};
		\end{scope}
		\begin{scope}[name prefix = t-]
			\node (1) at ($(m-1)+(0,2)$) {};
			\foreach \x in {2,...,\n}
			{
				\pgfmathtruncatemacro{\prec}{\x - 1}
				\node (\x) [right=of \prec] {};
			}
			\node (rdots)  [right=of \n] {\dots};
			\node (n)  [right=of rdots] {$n$};
			\node (ldots) [left=of 1] {\dots};
			\node (0) [left=of ldots] {1};
		\end{scope}
		
		\begin{scope}[name prefix = b-]
			\node (1) at ($(m-1)+(0,-3)$) {};
			\foreach \x in {2,...,\n}
			{
				\pgfmathtruncatemacro{\prec}{\x - 1}
				\node (\x) [right=of \prec] {};
			}
			\node (rdots)  [right=of \n] {\dots};
			\node (n)  [right=of rdots] {$n$};
			\node (ldots) [left=of 1] {\dots};
			\node (0) [left=of ldots] {1};
		\end{scope}
		
		\node at (t-1) {$j-1$};		
		\node at (b-1) {$j-1$};
		\node at (t-2) {$j$};
		\node at (t-3) {$j+1$};
		\node at (t-\n) {$i+2$};
		\node at (b-\nDeux) {$i$};
		\node at (b-\nUn) {$i+1$};
		\node at (b-\n) {$i+2$};
		
		\node at ($(m-0)+(-1,-1.5)$) {$d'$};
		\node at ($(m-0)+(-1.5,1)$) {$E_iE_{i-1} \dots E_j$};
		
		\path[{Circle[]}-{Circle[]}]
		(t-2) edge [bend right=90,looseness=2] (t-3)
		(b-\n) edge (t-\n)
		(b-n) edge (t-n);
		\path
		(m-\nDeux) edge [bend left=90,looseness=2,dotted,ForestGreen,thick] (m-\nUn)
		(b-0) edge[{Circle[]}-] ($(b-0)+(0,1)$)
		($(b-0)+(0,1)$) edge[dotted,thick] ($(b-0)+(0,1.5)$)
		($(m-0)+(0,0.1)$) edge[dotted,thick] ($(m-0)+(0,-0.5)$);
		
		\path[-{Circle[]}]
		(m-1) edge (t-1)
		(m-0) edge (t-0);
		
		\foreach \x in {2,...,\nTrois}
		{
			\pgfmathtruncatemacro{\suiv}{\x + 2}
			\draw[-{Circle[]}] (m-\x) edge[out=90,in=-90]  (t-\suiv);
		}
		
		\foreach \x in {1,...,\nTrois}
		{
			\path (b-\x) edge[{Circle[]}-] ($(b-\x)+(0,1)$)
			($(b-\x)+(0,1)$) edge[dotted,thick] ($(b-\x)+(0,1.5)$)
			($(m-\x)+(0,0.1)$) edge[dotted,thick] ($(m-\x)+(0,-0.5)$);
		}
		
		\path (b-\nUn) edge[{Circle[]}-,out=90,in=45] ($(m-\nDeux)+(0,-0.5)$)
		($(m-\nDeux)+(0,-0.5)$) edge[dotted,thick] ($(m-\nDeux)+(-0.3,-0.8)$)
		(b-\nUn) edge[dotted,ForestGreen,thick] ($(m-\nUn)+(0,0.1)$)
		($(m-\nDeux)+(0,0.1)$) edge[ForestGreen,dotted,thick] ($(m-\nDeux)+(0,-0.5)$)
		(b-\nDeux) edge[{Circle[]}-] ($(b-\nDeux)+(0,1)$)
		($(b-\nDeux)+(0,1)$) edge[dotted,thick] ($(b-\nDeux)+(0,1.5)$);
		
		\draw[dotted,red] (m-0) -- (m-n);
	\end{tikzpicture}
}

\newcommand{\figConfluence}{
	\begin{tikzpicture}[node distance=1cm,on grid]
		\node (0) at (0,0) {};
		\node (s) [below=of 0] {};
		\node (n) [above=of 0] {};
		\node (e) [right=of 0] {};
		\node (o) [left=of 0] {};
		
		\path[->>]
		(n) edge (e)
		(n) edge (o)
		(e) edge[dotted] (s)
		(o) edge[dotted] (s);
	\end{tikzpicture}
}

\newcommand{\figConfluenceLocale}{
	\begin{tikzpicture}[node distance=1cm,on grid]
		\node (0) at (0,0) {};
		\node (s) [below=of 0] {};
		\node (n) [above=of 0] {};
		\node (e) [right=of 0] {};
		\node (o) [left=of 0] {};
		
		\path
		(n) edge[->] (e)
		(n) edge[->] (o)
		(e) edge[dotted,->>] (s)
		(o) edge[dotted,->>] (s);
	\end{tikzpicture}
}

\newcommand{\figGrapheMot}{
	\begin{tikzpicture}[scale=0.8]
		\tikzset{
			grid_style/.style = {
				step=1cm,
				gray,
				very thin,
				dashed
			},
			cross/.style = {
				draw,
				blue,
				cross out,
				thick,
				minimum size=5pt,
				inner sep=0pt,
				outer sep=0pt,
				line width=1pt	
			}
		}
		
		\draw[grid_style] (0,0) grid (6,5);
		
		\path[->]
		(0,0) edge (7,0)
		edge (0,6);
		
		\foreach \x in {0,...,6}
		{
			\node at (\x,0) [below] {\x};
		}
		\node at (0,1) [left] {$\delta$};
		\node at (0,2) [left] {$e_1$};
		\node at (0,3) [left] {$e_2$};
		\node at (0,4) [left] {$e_3$};
		
		\node[cross] at (2,4) {};
		\node[cross] at (3,1) {};
		\node[cross] at (1,2) {};
		\node[cross] at (4,2) {};
	\end{tikzpicture} 
}

\newcommand{\figGrapheFNJ}{
	\begin{tikzpicture}[scale=0.8]
		\tikzset{
			grid_style/.style = {
				step=1cm,
				gray,
				very thin,
				dashed
			},
			cross/.style = {
				draw,
				blue,
				cross out,
				thick,
				minimum size=5pt,
				inner sep=0pt,
				outer sep=0pt,
				line width=1pt	
			}
		}
		
		\draw[grid_style] (0,0) grid (9,6);
		
		\path[->]
		(0,0) edge (10,0)
		edge (0,7);
		
		\foreach \x in {0,...,9}
		{
			\node at (\x,0) [below] {\x};
		}
		\node at (0,1) [left] {$\delta$};
		\node at (0,2) [left] {$e_1$};
		\node at (0,3) [left] {$e_2$};
		\node at (0,4) [left] {$e_3$};
		\node at (0,5) [left] {$e_4$};
		\node at (0,6) [left] {$e_5$};
		
		\node[cross] at (1,1) {};
		\node[cross] at (2,1) {};
		\node[cross] at (3,3) {};
		\node[cross] at (4,2) {};
		\node[cross] at (5,5) {};
		\node[cross] at (6,4) {};
		\node[cross] at (7,3) {};
		\node[cross] at (8,6) {};
		\node[cross] at (9,5) {};
		
		\draw[ForestGreen,dotted,thick] (4,2) -- (7,3) -- (9,5);
		\draw[ForestGreen,dotted,thick] (3,3) -- (5,5) -- (8,6);
		\path[blue,dotted,thick] 
		(4,2) edge (3,3)
		(5,5) edge (7,3)
		(8,6) edge (9,5);
	\end{tikzpicture}
}

\newcommand{\figQuatreDiagrammeOriente}{
	\begin{tikzpicture}[node distance=0.5cm]
		\begin{scope}[name prefix = t]
			\node (1) at (0,2) {1};
			\node (2) [right=of 1] {2};
			\node (3) [right=of 2] {3};
			\node (4) [right=of 3] {4};
		\end{scope}
		\begin{scope}[name prefix = b]
			\node (1) at (0,0) {1};
			\node (2) [right=of 1] {2};
			\node (3) [right=of 2] {3};
			\node (4) [right=of 3] {4};
		\end{scope}
		\path   (t1) edge[>->] (b1)
		(t2) edge [>->,bend right=90,looseness=2] (t3)
		(t4) edge [<-<,out=-90,in=90] (b2)
		(b3) edge [<-<,bend left=90,looseness=2] (b4);
	\end{tikzpicture}
}

\newcommand{\figRegleUn}{
	\begin{tikzpicture}[scale=0.7]
		\tikzset{
			grid_style/.style = {
				step=1cm,
				gray,
				very thin,
				dashed
			},
			cross/.style = {
				draw,
				blue,
				cross out,
				thick,
				minimum size=5pt,
				inner sep=0pt,
				outer sep=0pt,
				line width=1pt	
			}
		}
		
		\draw[grid_style] (0,0) grid (6,5);
		
		\path[->]
		(0,0) edge (7,0)
		edge (0,6);
		
		\foreach \x in {0,...,6}
		{
			\node at (\x,0) [below] {\x};
		}
		\node at (0,1) [left] {$\delta$};
		\node at (0,2) [left] {$e_1$};
		\node at (0,3) [left] {\dots};
		\node at (0,4) [left] {$e_i$};
		
		\node[cross] (1) at (2,4) {};
		\node[cross] (2) at (3,1) {};
		
		\draw[dotted,blue] (1) -- (2);
		\node at (7,3) {\huge{$\to$}};
		
		\node (o) at (9,0) {};
		\draw[grid_style] (o) grid ($(o)+(6,5)$);
		
		\path[->]
		(o) edge ($(o)+(7,0)$)
		edge ($(o)+(0,6)$);
		
		\foreach \x in {0,...,6}
		{
			\node at ($(o)+(\x,0)$) [below] {\x};
		}
		\node at ($(o)+(0,1)$) [left] {$\delta$};
		\node at ($(o)+(0,2)$) [left] {$e_1$};
		\node at ($(o)+(0,3)$) [left] {\dots};
		\node at ($(o)+(0,4)$) [left] {$e_i$};
		
		\node[cross] at ($(o)+(3,4)$) {};
		\node[cross] at ($(o)+(2,1)$) {};
	\end{tikzpicture} 
}

\newcommand{\figRegleDeux}{
	\begin{tikzpicture}[scale=0.7]
		\tikzset{
			grid_style/.style = {
				step=1cm,
				gray,
				very thin,
				dashed
			},
			cross/.style = {
				draw,
				blue,
				cross out,
				thick,
				minimum size=5pt,
				inner sep=0pt,
				outer sep=0pt,
				line width=1pt	
			}
		}
		
		\draw[grid_style] (0,0) grid (6,5);
		
		\path[->]
		(0,0) edge (7,0)
		edge (0,6);
		
		\foreach \x in {0,...,6}
		{
			\node at (\x,0) [below] {\x};
		}
		\node at (0,1) [left] {$\delta$};
		\node at (0,2) [left] {$e_1$};
		\node at (0,3) [left] {\dots};
		\node at (0,4) [left] {$e_i$};
		
		\node[cross] (1) at (2,4) {};
		\node[cross] (2) at (3,4) {};
		
		\draw[dotted,blue] (1) -- (2);
		\node at (7,3) {\huge{$\to$}};
		
		\node (o) at (9,0) {};
		\draw[grid_style] (o) grid ($(o)+(6,5)$);
		
		\path[->]
		(o) edge ($(o)+(7,0)$)
		edge ($(o)+(0,6)$);
		
		\foreach \x in {0,...,6}
		{
			\node at ($(o)+(\x,0)$) [below] {\x};
		}
		\node at ($(o)+(0,1)$) [left] {$\delta$};
		\node at ($(o)+(0,2)$) [left] {$e_1$};
		\node at ($(o)+(0,3)$) [left] {\dots};
		\node at ($(o)+(0,4)$) [left] {$e_i$};
		
		\node[cross] at ($(o)+(3,4)$) {};
		\node[cross] at ($(o)+(2,1)$) {};
	\end{tikzpicture} 
}

\newcommand{\figRegleTroisPlus}{
	\begin{tikzpicture}[scale=0.7]
		\tikzset{
			grid_style/.style = {
				step=1cm,
				gray,
				very thin,
				dashed
			},
			cross/.style = {
				draw,
				blue,
				cross out,
				thick,
				minimum size=5pt,
				inner sep=0pt,
				outer sep=0pt,
				line width=1pt	
			}
		}
		
		\draw[grid_style] (0,0) grid (6,5);
		
		\path[->]
		(0,0) edge (7,0)
		edge (0,6);
		
		\foreach \x in {0,...,6}
		{
			\node at (\x,0) [below] {\x};
		}
		\node at (0,1) [left] {\dots};
		\node at (0,2) [left] {\dots};
		\node at (0,3) [left] {$e_i$};
		\node at (0,4) [left] {$e_{i+1}$};
		
		\node[cross] (1) at (2,3) {};
		\node[cross] (2) at (3,4) {};
		\node[cross] (3) at (4,3) {};
		
		\draw[dotted,blue] (1) -- (2) -- (3);
		\node at (7,3) {\huge{$\to$}};
		
		\node (o) at (9,0) {};
		\draw[grid_style] (o) grid ($(o)+(6,5)$);
		
		\path[->]
		(o) edge ($(o)+(7,0)$)
		edge ($(o)+(0,6)$);
		
		\foreach \x in {0,...,6}
		{
			\node at ($(o)+(\x,0)$) [below] {\x};
		}
		\node at ($(o)+(0,1)$) [left] {\dots};
		\node at ($(o)+(0,2)$) [left] {\dots};
		\node at ($(o)+(0,3)$) [left] {$e_i$};
		\node at ($(o)+(0,4)$) [left] {$e_{i+1}$};
		
		\node[cross] at ($(o)+(2,3)$) {};
	\end{tikzpicture}
}

\newcommand{\figRegleTroisMoins}{
	\begin{tikzpicture}[scale=0.7]
		\tikzset{
			grid_style/.style = {
				step=1cm,
				gray,
				very thin,
				dashed
			},
			cross/.style = {
				draw,
				blue,
				cross out,
				thick,
				minimum size=5pt,
				inner sep=0pt,
				outer sep=0pt,
				line width=1pt	
			}
		}
		
		\draw[grid_style] (0,0) grid (6,5);
		
		\path[->]
		(0,0) edge (7,0)
		edge (0,6);
		
		\foreach \x in {0,...,6}
		{
			\node at (\x,0) [below] {\x};
		}
		\node at (0,1) [left] {\dots};
		\node at (0,2) [left] {\dots};
		\node at (0,3) [left] {$e_{i-1}$};
		\node at (0,4) [left] {$e_i$};
		
		\node[cross] (1) at (2,4) {};
		\node[cross] (2) at (3,3) {};
		\node[cross] (3) at (4,4) {};
		
		\draw[dotted,blue] (1) -- (2) -- (3);
		\node at (7,3) {\huge{$\to$}};
		
		\node (o) at (9,0) {};
		\draw[grid_style] (o) grid ($(o)+(6,5)$);
		
		\path[->]
		(o) edge ($(o)+(7,0)$)
		edge ($(o)+(0,6)$);
		
		\foreach \x in {0,...,6}
		{
			\node at ($(o)+(\x,0)$) [below] {\x};
		}
		\node at ($(o)+(0,1)$) [left] {\dots};
		\node at ($(o)+(0,2)$) [left] {\dots};
		\node at ($(o)+(0,3)$) [left] {$e_{i-1}$};
		\node at ($(o)+(0,4)$) [left] {$e_i$};
		
		\node[cross] at ($(o)+(2,4)$) {};
	\end{tikzpicture}
}

\newcommand{\figRegleQuatre}{
	\begin{tikzpicture}[scale=0.7]
		\tikzset{
			grid_style/.style = {
				step=1cm,
				gray,
				very thin,
				dashed
			},
			cross/.style = {
				draw,
				blue,
				cross out,
				thick,
				minimum size=5pt,
				inner sep=0pt,
				outer sep=0pt,
				line width=1pt	
			}
		}
		
		\draw[grid_style] (0,0) grid (6,5);
		
		\path[->]
		(0,0) edge (7,0)
		edge (0,6);
		
		\foreach \x in {0,...,6}
		{
			\node at (\x,0) [below] {\x};
		}
		\node at (0,1) [left] {\dots};
		\node at (0,2) [left] {$e_j$};
		\node at (0,3) [left] {\dots};
		\node at (0,4) [left] {$e_i$};
		
		\node[cross] (1) at (2,4) {};
		\node[cross] (2) at (3,2) {};
		
		\draw[dotted,blue] (1) -- (2);
		\node at (7,3) {\huge{$\to$}};
		
		\node (o) at (9,0) {};
		\draw[grid_style] (o) grid ($(o)+(6,5)$);
		
		\path[->]
		(o) edge ($(o)+(7,0)$)
		edge ($(o)+(0,6)$);
		
		\foreach \x in {0,...,6}
		{
			\node at ($(o)+(\x,0)$) [below] {\x};
		}
		\node at ($(o)+(0,1)$) [left] {\dots};
		\node at ($(o)+(0,2)$) [left] {$e_j$};
		\node at ($(o)+(0,3)$) [left] {\dots};
		\node at ($(o)+(0,4)$) [left] {$e_i$};
		
		\node[cross] at ($(o)+(2,2)$) {};
		\node[cross] at ($(o)+(3,4)$) {};
		
		\path[<->,thick,ForestGreen] (4,2) edge node[above right] {$\mathbf{\geq 2}$} (4,4);
	\end{tikzpicture} 
}

\newcommand{\figRegleCinq}{
	\begin{tikzpicture}[scale=0.7]
		\tikzset{
			grid_style/.style = {
				step=1cm,
				gray,
				very thin,
				dashed
			},
			cross/.style = {
				draw,
				blue,
				cross out,
				thick,
				minimum size=5pt,
				inner sep=0pt,
				outer sep=0pt,
				line width=1pt	
			}
		}
		
		\draw[grid_style] (0,0) grid (6,5);
		
		\path[->]
		(0,0) edge (7,0)
		edge (0,6);
		
		\foreach \x in {0,...,6}
		{
			\node at (\x,0) [below] {\x};
		}
		\node at (0,1) [left] {$e_{i-k}$};
		\node at (0,2) [left] {\dots};
		\node at (0,3) [left] {$e_{i-2}$};
		\node at (0,4) [left] {$e_{i-1}$};
		\node at (0,5) [left] {$e_i$};
		
		\node[cross] (0) at (1,5) {};
		\node[cross] (1) at (2,4) {};
		\node[cross] (2) at (3,3) {};
		\node[cross] (3) at (4,2) {};
		\node[cross] (4) at (5,1) {};
		\node[cross] (5) at (6,5) {};
		
		\draw[dotted,blue] (0) -- (1) -- (2) -- (3) -- (4) -- (5);
		\node at (7,3) {\huge{$\to$}};
		
		\node (o) at (9,0) {};
		\draw[grid_style] (o) grid ($(o)+(6,5)$);
		
		\path[->]
		(o) edge ($(o)+(7,0)$)
		edge ($(o)+(0,6)$);
		
		\foreach \x in {0,...,6}
		{
			\node at ($(o)+(\x,0)$) [below] {\x};
		}
		\node at ($(o)+(0,1)$) [left] {$e_{i-k}$};
		\node at ($(o)+(0,2)$) [left] {\dots};
		\node at ($(o)+(0,3)$) [left] {$e_{i-2}$};
		\node at ($(o)+(0,4)$) [left] {$e_{i-1}$};
		\node at ($(o)+(0,5)$) [left] {$e_i$};
		
		\node[cross] at ($(o)+(1,3)$) {};
		\node[cross] at ($(o)+(2,2)$) {};
		\node[cross] at ($(o)+(3,1)$) {};
		\node[cross] at ($(o)+(4,5)$) {};
	\end{tikzpicture}  
}

\newcommand{\figRegleSix}{
	\begin{tikzpicture}[scale=0.7]
		\tikzset{
			grid_style/.style = {
				step=1cm,
				gray,
				very thin,
				dashed
			},
			cross/.style = {
				draw,
				blue,
				cross out,
				thick,
				minimum size=5pt,
				inner sep=0pt,
				outer sep=0pt,
				line width=1pt	
			}
		}
		
		\draw[grid_style] (0,0) grid (6,5);
		
		\path[->]
		(0,0) edge (7,0)
		edge (0,6);
		
		\foreach \x in {0,...,6}
		{
			\node at (\x,0) [below] {\x};
		}
		\node at (0,1) [left] {$e_i$};
		\node at (0,2) [left] {$e_{i+1}$};
		\node at (0,3) [left] {$e_{i+2}$};
		\node at (0,4) [left] {\dots};
		\node at (0,5) [left] {$e_{i+k}$};
		
		\node[cross] (0) at (1,1) {};
		\node[cross] (1) at (2,5) {};
		\node[cross] (2) at (3,4) {};
		\node[cross] (3) at (4,3) {};
		\node[cross] (4) at (5,2) {};
		\node[cross] (5) at (6,1) {};
		
		\draw[dotted,blue] (0) -- (1) -- (2) -- (3) -- (4) -- (5);
		\node at (7,3) {\huge{$\to$}};
		
		\node (o) at (9,0) {};
		\draw[grid_style] (o) grid ($(o)+(6,5)$);
		
		\path[->]
		(o) edge ($(o)+(7,0)$)
		edge ($(o)+(0,6)$);
		
		\foreach \x in {0,...,6}
		{
			\node at ($(o)+(\x,0)$) [below] {\x};
		}
		\node at ($(o)+(0,1)$) [left] {$e_i$};
		\node at ($(o)+(0,2)$) [left] {$e_{i+1}$};
		\node at ($(o)+(0,3)$) [left] {$e_{i+2}$};
		\node at ($(o)+(0,4)$) [left] {\dots};
		\node at ($(o)+(0,5)$) [left] {$e_{i+k}$};
		
		\node[cross] at ($(o)+(1,1)$) {};
		\node[cross] at ($(o)+(2,5)$) {};
		\node[cross] at ($(o)+(3,4)$) {};
		\node[cross] at ($(o)+(4,3)$) {};
	\end{tikzpicture} 
}
\newcommand{\figW}{
	\begin{tikzpicture}[scale=0.7]
		\tikzset{
			grid_style/.style = {step=1cm,gray,very thin,dashed},
			cross/.style = {draw,blue,cross out,thick,minimum size=5pt,inner sep=0pt,outer sep=0pt,line width=1pt}
		}
		\draw[grid_style] (0,0) grid (16,6);
		
		\path[->]
		(0,0) edge (17,0)
		edge (0,7);
		
		\foreach \x in {0,...,16}
		{
			\node at (\x,0) [below] {\x};
		}
		\node at (0,1) [left] {$\delta$};
		\node at (0,2) [left] {$e_1$};
		\node at (0,3) [left] {$e_2$};
		\node at (0,4) [left] {$e_3$};
		\node at (0,5) [left] {$e_4$};
		\node at (0,6) [left] {$e_5$};
		
		\path[red,dashed,line width = 1.5pt]
		(0,5) edge (16,5)
		(3,0) edge (3,5)
		(10,0) edge (10,5)
		(12,0) edge (12,5);
		
		\node[cross] at (1,3) {};
		\node[cross] at (2,2) {};
		\node[cross] at (3,5) {};
		\node[cross] at (4,4) {};
		\node[cross] at (5,3) {};
		\node[cross] at (6,4) {};
		\node[cross] at (7,1) {};
		\node[cross] at (8,4) {};
		\node[cross] at (9,4) {};
		\node[cross] at (10,5) {};
		\node[cross] at (11,2) {};
		\node[cross] at (12,5) {};
		\node[cross] at (13,4) {};
		\node[cross] at (14,2) {};
		\node[cross] at (15,3) {};
		
		\path 
		(0,-0.5) -- node[below,ForestGreen] {{$\bm{v_0}$}} (3,-0.5) node[below,red] {{$\bm{e_m}$}}
		(3,-0.5) -- node[below,ForestGreen] {$\bm{v_1}$} (10,-0.5) node[below,red] {{$\bm{e_m}$}}
		(10,-0.5) -- node[below,ForestGreen] {$\bm{v_2}$} (12,-0.5) node[below,red] {{$\bm{e_m}$}}
		(12,-0.5) -- node[below,ForestGreen] {$\bm{v_3}$} (15,-0.5);
		
	\end{tikzpicture}
}

\newcommand{\figVprime}{
	\begin{tikzpicture}[scale=0.7]
		\tikzset{
			grid_style/.style = {step=1cm,gray,very thin,dashed},
			cross/.style = {draw,blue,cross out,thick,minimum size=5pt,inner sep=0pt,outer sep=0pt,line width=1pt}
		}
		\draw[grid_style] (0,0) grid (16,6);
		
		\path[->]
		(0,0) edge (17,0)
		edge (0,7);
		
		\foreach \x in {0,...,16}
		{
			\node at (\x,0) [below] {\x};
		}
		\node at (0,1) [left] {$\delta$};
		\node at (0,2) [left] {$e_1$};
		\node at (0,3) [left] {$e_2$};
		\node at (0,4) [left] {$e_3$};
		\node at (0,5) [left] {$e_4$};
		\node at (0,6) [left] {$e_5$};
		
		\path[red,dashed,line width = 1.5pt]
		(0,5) edge (16,5)
		(3,0) edge (3,5)
		(10,0) edge (10,5)
		(12,0) edge (12,5);
		
		\node[cross] at (1,3) {};
		\node[cross] at (2,2) {};
		\node[cross] at (3,5) {};
		\node[cross] (1) at (4,1) {};
		\node[cross] (2) at (5,1) {};
		\node[cross] (3) at (6,3) {};
		\node[cross] (4) at (7,2) {};
		\node[cross] (5) at (8,4) {};
		\node[cross] (6) at (9,3) {};
		\node[cross] at (10,5) {};
		\node[cross] at (11,2) {};
		\node[cross] at (12,5) {};
		\node[cross] at (13,4) {};
		\node[cross] at (14,2) {};
		\node[cross] at (15,3) {};
		
		\path (3) edge[dotted,blue,thick] (4)
		(5) edge[dotted,blue,thick] (6)
		(4) edge[dotted,ForestGreen,thick] (6)
		(3) edge[dotted,ForestGreen,thick] (5);
		\path 
		(0,-0.5) -- node[below,ForestGreen] {{$\bm{v_0}$}} (3,-0.5) node[below,red] {{$\bm{e_m}$}}
		(3,-0.5) -- node[below,ForestGreen] {$\bm{v_1'}$} (10,-0.5) node[below,red] {{$\bm{e_m}$}}
		(10,-0.5) -- node[below,ForestGreen] {$\bm{v_2}$} (12,-0.5) node[below,red] {{$\bm{e_m}$}}
		(12,-0.5) -- node[below,ForestGreen] {$\bm{v_3}$} (15,-0.5);
		
	\end{tikzpicture}
}

\newcommand{\figSdeuxVseconde}{
	\begin{tikzpicture}[scale=0.7]
		\tikzset{
			grid_style/.style = {step=1cm,gray,very thin,dashed},
			cross/.style = {draw,blue,cross out,thick,minimum size=5pt,inner sep=0pt,outer sep=0pt,line width=1pt}
		}
		\draw[grid_style] (0,0) grid (16,6);
		
		\path[->]
		(0,0) edge (17,0)
		edge (0,7);
		
		\foreach \x in {0,...,16}
		{
			\node at (\x,0) [below] {\x};
		}
		\node at (0,1) [left] {$\delta$};
		\node at (0,2) [left] {$e_1$};
		\node at (0,3) [left] {$e_2$};
		\node at (0,4) [left] {$e_3$};
		\node at (0,5) [left] {$e_4$};
		\node at (0,6) [left] {$e_5$};
		
		\path[red,dashed,line width = 1.5pt]
		(0,5) edge (16,5)
		(7,0) edge (7,5)
		(10,0) edge (10,5)
		(12,0) edge (12,5);
		
		\node[cross] at (1,3) {};
		\node[cross] at (2,2) {};
		\node[cross] at (3,1) {};
		\node[cross] at (4,1) {};
		\node[cross] at (5,3) {};
		\node[cross] at (6,2) {};
		\node[cross] (1) at (7,5) {};
		\node[cross] (2) at (8,4) {};
		\node[cross] (3) at (9,3) {};
		\node[cross] (4) at (10,5) {};
		\node[cross] at (11,2) {};
		\node[cross] at (12,5) {};
		\node[cross] at (13,4) {};
		\node[cross] at (14,2) {};
		\node[cross] at (15,3) {};
		
		\draw[dotted,blue,thick] (1)--(3)--(4);
		\path 
		(0,-0.5) -- node[below,ForestGreen] {{$\bm{v_0'}$}} (7,-0.5) node[below,red] {{$\bm{e_m}$}}
		(7,-0.5) -- node[below,ForestGreen] {$\bm{v_1''}$} (10,-0.5) node[below,red] {{$\bm{e_m}$}}
		(10,-0.5) -- node[below,ForestGreen] {$\bm{v_2}$} (12,-0.5) node[below,red] {{$\bm{e_m}$}}
		(12,-0.5) -- node[below,ForestGreen] {$\bm{v_3}$} (15,-0.5);
		
	\end{tikzpicture}
}

\newcommand{\figSunVseconde}{
	\begin{tikzpicture}[scale=0.7]
		\tikzset{
			grid_style/.style = {step=1cm,gray,very thin,dashed},
			cross/.style = {draw,blue,cross out,thick,minimum size=5pt,inner sep=0pt,outer sep=0pt,line width=1pt}
		}
		\draw[grid_style] (0,0) grid (13,7);
		
		\path[->]
		(0,0) edge (14,0)
		edge (0,8);
		
		\foreach \x in {0,...,13}
		{
			\node at (\x,0) [below] {\x};
		}
		\node at (0,1) [left] {$\delta$};
		\node at (0,2) [left] {$e_1$};
		\node at (0,3) [left] {$e_2$};
		\node at (0,4) [left] {$e_3$};
		\node at (0,5) [left] {$e_4$};
		\node at (0,6) [left] {$e_5$};
		
		\path[red,dashed,line width = 1.5pt]
		(0,6) edge (13,6)
		(8,0) edge (8,6);
		
		\node[cross] (1) at (1,1) {};
		\node[cross] (2) at (2,1) {};
		\node[cross] (3) at (3,3) {};
		\node[cross] (4) at (4,2) {};
		\node[cross] (5) at (5,5) {};
		\node[cross] (6) at (6,4) {};
		\node[cross] (7) at (7,3) {};
		\node[cross] (8) at (8,6) {};
		\node[cross] (9) at (9,5) {};
		\node[cross] (10)at (10,4) {};
		\node[cross] (11)at (11,3) {};
		\node[cross,violet] (12)at (12,2) {};
		
		\path[dotted,blue,thick] 
		(3) edge (4)
		(5) edge (7)
		(8) edge (12);
		\path 
		(0,-0.5) -- node[below,ForestGreen] {{$\bm{v_0''}$}} (8,-0.5) node[below,red] {{$\bm{e_m}$}}
		(8,-0.5) -- node[below,ForestGreen] {$\bm{v_1''}$} (12,-0.5) ;
		
		\path[ForestGreen,dotted,thick] 
		(3) edge (5)
		(4) edge (7)
		(5) edge (8);
		\path[dotted,thick]
		(7) edge (12);
	\end{tikzpicture}
}

\newcommand{\figSunU}{
	\begin{tikzpicture}[scale=0.7]
		\tikzset{
			grid_style/.style = {step=1cm,gray,very thin,dashed},
			cross/.style = {draw,blue,cross out,thick,minimum size=5pt,inner sep=0pt,outer sep=0pt,line width=1pt}
		}
		\draw[grid_style] (0,0) grid (11,7);
		
		\path[->] (0,0) 
		edge (12,0)
		edge (0,8);
		
		\foreach \x in {0,...,13}
		{
			\node at (\x,0) [below] {\x};
		}
		\node at (0,1) [left] {$\delta$};
		\node at (0,2) [left] {$e_1$};
		\node at (0,3) [left] {$e_2$};
		\node at (0,4) [left] {$e_3$};
		\node at (0,5) [left] {$e_4$};
		\node at (0,6) [left] {$e_5$};
		
		\node[cross] (1) at (1,1) {};
		\node[cross] (2) at (2,1) {};
		\node[cross] (3) at (3,3) {};
		\node[cross] (4) at (4,2) {};
		\node[cross] (5) at (5,5) {};
		\node[cross] (6) at (6,4) {};
		\node[cross] (7) at (7,3) {};
		\node[cross,violet] (8) at (8,2) {};
		\node[cross] (9) at (9,6) {};
		\node[cross] (10)at (10,5) {};
		
		\path[dotted,blue,thick] 
		(3) edge (4)
		(5) edge (8)
		(9) edge (10);
		\path 
		(0,-0.5) -- node[below,ForestGreen] {{$\bm{v_0''}$}} (8,-0.5) node[below,violet] {{$\bm{u_1}$}}
		(8,-0.5) -- node[below] {$\bm{u_0}$} (11,-0.5) ;
		
		\path[ForestGreen,dotted,thick] 
		(3) edge (5)
		(8) edge (10)
		(5) edge (9);
		\path[dashed,line width = 1.5pt]
		(4) edge (8);
	\end{tikzpicture}
}
 \newcommand{\dist}{5}

\newcommand{\figCap}[1][1]{%
	\vcenter{\hbox{%
			\begin{tikzpicture}[node distance=\dist mm, xscale=#1]
				\node (0) {};
				\node (1) [right=of 0] {};
				\path[{Circle[scale=0.5]}-{Circle[scale=0.5]}] (0) edge[bend left=90, looseness=2] (1);
			\end{tikzpicture}%
	}}%
}
\newcommand{\figCup}{
	\figCap[-1]
}

\newcommand{\figId}[1][\dist]{
	\vcenter{\hbox{
			\begin{tikzpicture}[node distance = #1 mm]
				\node (0) {};
				\node (1) [above=of 0] {};
				\path[{Circle[scale=0.5]}-{Circle[scale=0.5]}] (0) edge (1);
			\end{tikzpicture}
	}}
}

\newcommand{\figZigzag}[1][1]{
	\pgfmathsetmacro{\halfdist}{\dist/2}
	\pgfmathsetmacro{\mhalfdist}{-\halfdist}
	\vcenter{\hbox{
			\begin{tikzpicture}[node distance = \halfdist mm,xscale= #1,x=1mm, y=1mm]
				\node (2) at (0,0) {};
				\node (1) at (\mhalfdist,0) {};
				\node (0) at (\mhalfdist,\mhalfdist) {};
				\node (3) at (\halfdist,0) {};
				\node (4) at (\halfdist,\halfdist) {};
				\path
				(0.center) edge[{Circle[scale=0.5]}-] (1.center)
				(1.center) edge[bend left =90,looseness=2] (2.center)
				(2.center) edge[bend right =90,looseness=2] (3.center)
				(3.center) edge[-{Circle[scale=0.5]}] (4.center)
				;
			\end{tikzpicture}
	}}
}

\newcommand{\figBubble}{
	\vcenter{\hbox{
			\begin{tikzpicture}
				\pgfmathsetmacro{\halfdist}{\dist/2}
				\draw (0,0) circle (\halfdist mm);
			\end{tikzpicture}
	}}
}

\newcommand{\figEchangeA}{
	\vcenter{\hbox{
			\begin{tikzpicture}[node distance = 10 mm and 5mm]
				\node (1-1) {$n$};
				\node (1-2) [above=of 1-1] {$n$};
				\node (1-3) [above=of 1-2] {$p$};
				\node (2-1) [right=of 1-1] {$m$};
				\node (2-2) [above=of 2-1] {$m$};
				\node (2-3) [above=of 2-2] {$q$};
				\path
				(1-1) edge[->] node[right] {$f$} (1-2)
				(1-2) edge[->] node[right,darkgray] {$1_p$} (1-3)
				(2-1) edge[->] node[right,darkgray] {$1_m$} (2-2)
				(2-2) edge[->] node[right] {$g$} (2-3)
				;
			\end{tikzpicture}
	}}
}

\newcommand{\figEchangeB}{
	\vcenter{\hbox{
			\begin{tikzpicture}[node distance = 10 mm and 5mm]
				\node (1-1) {$n$};
				\node (1-2) [above=of 1-1] {$n$};
				\node (1-3) [above=of 1-2] {$p$};
				\node (2-1) [right=of 1-1] {$m$};
				\node (2-2) [above=of 2-1] {$q$};
				\node (2-3) [above=of 2-2] {$q$};
				\path
				(1-1) edge[->] node[right,darkgray] {$1_n$} (1-2)
				(1-2) edge[->] node[right] {$f$} (1-3)
				(2-1) edge[->] node[right] {$g$} (2-2)
				(2-2) edge[->] node[right,darkgray] {$1_q$} (2-3)
				;
			\end{tikzpicture}
	}}
}

\newcommand{\gdist}{10}
\newcommand{\figCapP}[1][1]{%
	\vcenter{\hbox{%
			\begin{tikzpicture}[node distance=\gdist mm, xscale=#1]
				\node (0) {};
				\node (1) [right=of 0] {};
				\path[{Straight Barb[length=8pt,width=8pt]}-{Straight Barb[length=8pt,width=8pt,reversed]},
				shorten <=8pt]
				(0) edge[bend left=90, looseness=2.5] (1);
			\end{tikzpicture}%
	}}%
}
\newcommand{\figCupP}{
	\figCapP[-1]
}

\newcommand{\figCapM}[1][1]{%
	\vcenter{\hbox{%
			\begin{tikzpicture}[node distance=\gdist mm,xscale=#1]
				\node (0) {};
				\node (1) [right=of 0] {};
				\path[{Straight Barb[length=8pt,width=8pt,reversed]}-{Straight Barb[length=8pt,width=8pt]},
				shorten >=8pt]
				(0) edge[bend left=90, looseness=2.5] (1);
			\end{tikzpicture}%
	}}%
}
\newcommand{\figCupM}{
	\figCapM[-1]
}

\newcommand{\figZigzagUp}[1][1]{
	\pgfmathsetmacro{\halfdist}{\gdist/2}
	\pgfmathsetmacro{\mhalfdist}{-\halfdist}
	\pgfmathsetmacro{\mgdist}{-\gdist}
	\vcenter{\hbox{
			\begin{tikzpicture}[node distance = \halfdist mm,xscale= #1,x=1mm, y=1mm]
				\node (2) at (0,0) {};
				\node (1) at (\mhalfdist,0) {};
				\node (0) at (\mhalfdist,\mgdist) {};
				\node (3) at (\halfdist,0) {};
				\node (4) at (\halfdist,\gdist) {};
				\path
				(0.center) edge[{Straight Barb[length=8pt,width=8pt,reversed]}-] (1.center)
				(1.center) edge[bend left =90,looseness=2] (2.center)
				(2.center) edge[bend right =90,looseness=2] (3.center)
				(3.center) edge[-{Straight Barb[length=8pt,width=8pt]}] (4.center)
				;
			\end{tikzpicture}
	}}
}

\newcommand{\figZigzagDown}[1][1]{
	\pgfmathsetmacro{\halfdist}{\gdist/2}
	\pgfmathsetmacro{\mhalfdist}{-\halfdist}
	\pgfmathsetmacro{\mgdist}{-\gdist}
	\vcenter{\hbox{
			\begin{tikzpicture}[node distance = \halfdist mm,xscale= #1,x=1mm, y=1mm]
				\node (2) at (0,0) {};
				\node (1) at (\mhalfdist,0) {};
				\node (0) at (\mhalfdist,\mgdist) {};
				\node (3) at (\halfdist,0) {};
				\node (4) at (\halfdist,\gdist) {};
				\path
				(0.center) edge[{Straight Barb[length=8pt,width=8pt]}-] (1.center)
				(1.center) edge[bend left =90,looseness=2] (2.center)
				(2.center) edge[bend right =90,looseness=2] (3.center)
				(3.center) edge[-{Straight Barb[length=8pt,width=8pt,reversed]}] (4.center)
				;
			\end{tikzpicture}
	}}
}

\newcommand{\figIdUp}{
	\vcenter{\hbox{
			\begin{tikzpicture}[x=1mm,y=1mm]
				\node (0) at (0,0) {};
				\node (1) at (0,24) {};
				\path[{Straight Barb[length=8pt,width=8pt,reversed]}-{Straight Barb[length=8pt,width=8pt]}] (0) edge (1);
			\end{tikzpicture}
	}}
}

\newcommand{\figIdDown}{
	\vcenter{\hbox{
			\begin{tikzpicture}[x=1mm,y=1mm]
				\node (0) at (0,0) {};
				\node (1) at (0,24) {};
				\path[{Straight Barb[length=8pt,width=8pt]}-{Straight Barb[length=8pt,width=8pt,reversed]}] (0) edge (1);
			\end{tikzpicture}
	}}
}

\newcommand{\figBubbleTrigo}{
	\vcenter{\hbox{
			\begin{tikzpicture}
				\node at (5mm,0) {$\wedge$};
				\node at (-5mm,0) {$\vee$};
				\draw (0,0) circle (5 mm);
			\end{tikzpicture}
	}}
}

\newcommand{\figBubbleHor}{
	\vcenter{\hbox{
			\begin{tikzpicture}
				\node at (-5mm,0) {$\wedge$};
				\node at (5mm,0) {$\vee$};
				\draw (0,0) circle (5 mm);
			\end{tikzpicture}
	}}
}

\newcommand{\figConfluenceLocaleModulo}{
	\begin{tikzpicture}[node distance=1cm,on grid]
		\node (0) at (0,0) {};
		\node (g1) [below left=of 0] {};
		\node (g2) [below left=of g1] {};
		\node (g3) [below left=of g2] {};
		\node (g4) [below right=of g3] {};
		\node (g5) [below right=of g4] {};
		\node (g6) [below right=of g5] {};
		\node (d1) [below right=of 0] {};
		\node (d2) [below right=of d1] {};
		\node (d3) [below right=of d2] {};
		\node (d4) [below left=of d3] {};
		\node (d5) [below left=of d4] {};
		\node (d6) [below left=of d5] {};
		
		\path
		(0) edge[decorate, decoration={snake, amplitude=0.5mm, segment length=3mm},->] (g1)
		edge[decorate, decoration={snake, amplitude=0.5mm, segment length=3mm},->] (d1)
		(g1) edge[->] (g2)
		(d1) edge[->] (d2)
		(g2) edge[decorate, decoration={snake, amplitude=0.5mm, segment length=3mm},->] (g3)
		(d2) edge[decorate, decoration={snake, amplitude=0.5mm, segment length=3mm},->] (d3)
		(g3) edge[decorate, decoration={snake, amplitude=0.5mm, segment length=3mm},->,dotted,thick] (g4)
		(d3) edge[decorate, decoration={snake, amplitude=0.5mm, segment length=3mm},->,dotted,thick] (d4)
		(g4) edge[->>,dotted,thick] (g5)
		(d4) edge[->>,dotted,thick] (d5)
		(g5) edge[decorate, decoration={snake, amplitude=0.5mm, segment length=3mm},->,dotted,thick] (g6)
		(d5) edge[decorate, decoration={snake, amplitude=0.5mm, segment length=3mm},->,dotted,thick] (d6)
		;
	\end{tikzpicture}
}

\newcommand{\figConfluenceLocaleModuloBis}{
	\begin{tikzpicture}[node distance=1cm,on grid]
		\node (0) at (0,0) {};
		\node (g1) [below left=of 0] {};
		\node (g2) [below left=of g1] {};
		\node (g3) [below left=of g2] {};
		\node (g4) [below right=of g3] {};
		\node (g5) [below right=of g4] {};
		\node (g6) [below right=of g5] {};
		\node (d1) [below right=of 0] {};
		\node (d2) [below right=of d1] {};
		\node (d3) [below right=of d2] {};
		\node (d4) [below left=of d3] {};
		\node (d5) [below left=of d4] {};
		\node (d6) [below left=of d5] {};
		
		\path
		(g2) edge[decorate, decoration={snake, amplitude=0.5mm, segment length=3mm},<->] (d2)
		(g2) edge[->] (g3)
		(d2) edge[->] (d3)
		(g3) edge[decorate, decoration={snake, amplitude=0.5mm, segment length=3mm},->,dotted,thick] (g4)
		(d3) edge[decorate, decoration={snake, amplitude=0.5mm, segment length=3mm},->,dotted,thick] (d4)
		(g4) edge[->>,dotted,thick] (g5)
		(d4) edge[->>,dotted,thick] (d5)
		(g5) edge[decorate, decoration={snake, amplitude=0.5mm, segment length=3mm},->,dotted,thick] (g6)
		(d5) edge[decorate, decoration={snake, amplitude=0.5mm, segment length=3mm},->,dotted,thick] (d6)
		;
	\end{tikzpicture}
}

\newcommand{\figConfluenceCategorie}{
	\begin{tikzpicture}[node distance=3cm and 4cm,on grid]   
		\node (o) {
			\centering
			\begin{tikzpicture}[node distance=5mm]
				\node (0) at (0,0) {};
				\node (1) [above=of 0] {};
				\node (2) [right=of 1] {};
				\node (3) [right=of 2] {};
				\node (4) [above=of 3] {};
				\node (5) [right=of 4] {};
				\node (6) [below=of 5] {};
				\node (7) [below=of 6] {};
				
				\draw[ForestGreen,{Straight Barb[length=8pt,width=8pt,reversed]}-] (0.center) to (1.center) to[bend left=90,looseness=2] (2.center) to[bend right=90,looseness=2] (3.center) to (4.center);
				\draw[blue,-{Straight Barb[length=8pt,width=8pt]}] (4.center) to[bend left=90,looseness=2] (5.center) to (6.center) to[blue] (7.center);
			\end{tikzpicture}
		};
		
		\node[right=of o] (d) {
			\centering
			\begin{tikzpicture}[node distance=5mm]
				\node (0) at (0,0) {};
				\node (1) [above=of 0] {};
				\node (2) [above=of 1] {};
				\node (3) [right=of 2] {};
				\node (4) [below=of 3] {};
				\node (5) [right=of 4] {};
				\node (6) [right=of 5] {};
				\node (7) [below=of 6] {};
				
				\draw[ForestGreen,{Straight Barb[length=8pt,width=8pt,reversed]}-] (0.center) to (1.center) to (2.center) to[bend left=90,looseness=2] (3.center);
				\draw[blue,-{Straight Barb[length=8pt,width=8pt]}] (3.center) to (4.center) to[bend right=90,looseness=2] (5.center) to[bend left=90,looseness=2] (6.center) to[blue] (7.center);
			\end{tikzpicture}
		};
		
		\node[below=of o] (b) {
			\centering
			\begin{tikzpicture}[node distance=5mm]
				\node (0) at (0,0) {};
				\node (1) [above=of 0] {};
				\node (2) [above=of 1] {};
				\node (3) [right=of 2] {};
				\node (4) [below=of 3] {};
				\node (5) [below=of 4] {};
				
				\draw[ForestGreen,{Straight Barb[length=8pt,width=8pt,reversed]}-] (0.center) to (1.center) to (2.center);
				\draw[blue,-{Straight Barb[length=8pt,width=8pt]}] (2.center) to[bend left=90,looseness=2] (3.center) to (4.center) to (5.center);
			\end{tikzpicture}
		};
		
		\node[right=of b] (db) {
			\centering
			\begin{tikzpicture}[node distance=5mm]
				\node (0) at (0,0) {};
				\node (1) [above=of 0] {};
				\node (2) [above=of 1] {};
				\node (3) [right=of 2] {};
				\node (4) [below=of 3] {};
				\node (5) [below=of 4] {};
				
				\draw[ForestGreen,{Straight Barb[length=8pt,width=8pt,reversed]}-] (0.center) to (1.center) to (2.center) to[bend left=90,looseness=2] (3.center);
				\draw[blue,-{Straight Barb[length=8pt,width=8pt]}] (3.center) to (4.center) to (5.center);
			\end{tikzpicture}
		};
		\draw[decorate, decoration={snake, amplitude=0.5mm, segment length=3mm},<->] (o) -- (d) node[midway,above] {modulo};
		\draw[->] (d) -- (db) node[midway,right] {(1-right)};
		\draw[->] (o) -- (b) node[midway,left] {(1-left)};
		\draw[dotted] (b) -- (db) node[midway,above] {\huge{$=$}};
		
	\end{tikzpicture}
}

\usepackage[english]{babel}

\usepackage[T1]{fontenc}
\usepackage[utf8]{inputenc}
\usepackage[a4paper,top=2cm,bottom=2cm,left=2cm,right=2cm,marginparwidth=1.75cm]{geometry}

\usepackage[backend=biber,sorting=nyt,style=alphabetic]{biblatex}
\usepackage{graphicx}
\usepackage[colorlinks=true, allcolors=blue]{hyperref}
\usepackage{algorithm}
\usepackage[noend]{algpseudocode}
\usepackage{titlesec}
\usepackage[dvipsnames]{xcolor}
\usepackage{placeins}
\usepackage{subcaption}
\usepackage[framemethod=tikz]{mdframed}
\usepackage[all]{xy}
\usepackage{bm}
\usepackage{xcolor}
\titleformat{\section}[block]
{\normalfont\Large\bfseries\centering}
{\thesection}         
{1em}                 
{}                    
[                    
\vspace{-4ex}     
\titlerule          
\vspace{5ex}
\titlerule          
\vspace{1ex}
]

\usepackage{amsmath}
\usepackage{amsthm}
\usepackage{amssymb}

\usepackage{tikz}
\usetikzlibrary{decorations.pathreplacing}
\usetikzlibrary{positioning}
\usetikzlibrary{arrows.meta}
\usetikzlibrary{calc}
\usetikzlibrary {fit,shapes.geometric}
\usetikzlibrary{matrix}
\usetikzlibrary{shapes.misc}
\usetikzlibrary{decorations.pathmorphing}

\theoremstyle{plain}
\newmdtheoremenv[hidealllines=true,leftline=true,innerleftmargin=10pt,innerrightmargin=0pt,innertopmargin=0pt]{theorem}{Theorem}[section]
\newmdtheoremenv[hidealllines=true,leftline=true,innerleftmargin=10pt,innerrightmargin=0pt,innertopmargin=0pt]{corollary}{Corollary}[theorem]
\newmdtheoremenv[hidealllines=true,leftline=true,innerleftmargin=10pt,innerrightmargin=0pt,innertopmargin=0pt]{lemma}[theorem]{Lemma}
\newmdtheoremenv[hidealllines=true,leftline=true,innerleftmargin=10pt,innerrightmargin=10pt,innertopmargin=0pt]{proposition}[theorem]{Proposition}
\newmdtheoremenv[hidealllines=true,leftline=true,innerleftmargin=10pt,innerrightmargin=0pt,innertopmargin=0pt]{conjecture}[theorem]{Conjecture}

\theoremstyle{definition}
\newmdtheoremenv[hidealllines=true,leftline=true,innerleftmargin=10pt,innerrightmargin=0pt,innertopmargin=0pt]{definition}[theorem]{Definition}
\newtheorem*{notation}{Notation}
\newtheorem*{vocabulary}{Vocabulary}
\newtheorem{example}[theorem]{Example}
\newtheorem{remark}[theorem]{Remark}

\makeatletter
\renewcommand{\ALG@name}{\textsc{Algorithm}}
\makeatother
\numberwithin{algorithm}{section}

\newcommand{\cO}{\mathcal{O}}
\newcommand{\cF}{\mathcal{F}}
\newcommand{\N}{\mathbb{N}}
\newcommand{\cTL}{\mathcal{T\!L}}
\newcommand{\green}[1]{\textcolor{ForestGreen}{#1}}
\newcommand{\red}[1]{\textcolor{red}{#1}}

\begin{document}
	\begin{par}
		\begin{figure}[t]
			\centering
			\includegraphics[width=0.25\linewidth]{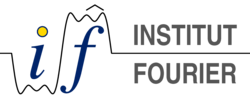}
			\hspace{0.15\linewidth}
			\includegraphics[width=0.25\linewidth]{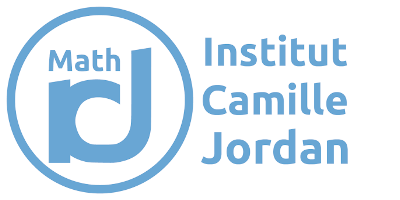}
		\end{figure}
		\begin{centering}
			{\Large \href{https://arxiv.org/abs/2508.19360}{\color{black}{Internship report}}\par}
			{\huge \textbf{Search for a basis of the Temperley-Lieb algebra, using rewriting systems}\par}
		\end{centering}
		\vspace{1cm}
		{Done by : Julien THIÉBAUT\par}
		{Mail : julien.thiebaut [at] etu.univ-grenoble-alpes [dot] fr\par}
		{\textsc{Bachelor's degree in mathematics}\par}
		{\textsc{Institut Joseph Fourier}\par}
		{\textsc{Université Grenoble-Alpes}\par}
		
		\vspace{0.5cm}
		{Tutor : Stéphane GAUSSENT\par}
		{\textsc{Teacher-researcher}\par}
		{\textsc{Institut Camille Jordan}\par}
		{\textsc{Université Jean Monnet}\par}
		\vspace{0.5cm}
		{\textit{From June 2nd to July 11th 2025.}\par}
	\end{par}
	
	\vspace{1cm} 
	
	\begin{figure}[ht]
		\centering
		\figPremierePage
		\label{fig:premPage}
	\end{figure}
	
	\vspace{1cm} 
	
	\begin{abstract}
		We begin by defining Temperley-Lieb algebra, in two different ways: as a presented algebra or as a diagrammatic algebra. Next, we look for a basis algorithmically, using rewriting theory. Finally, we introduce a generalization of the Temperley-Lieb algebra, which is an oriented version of the previous one. This pushes us to employ a more efficient tool, category theory, to use rewriting to easily obtain a basis for the algebra.
	\end{abstract}
	
	\newpage
	\tableofcontents
	
	\newpage
	\section*{Introduction}
	
	To complete my Bachelor's degree in mathematics in Grenoble, I decided to do a research internship in Saint-Étienne, my hometown. I'd like to thank Stéphane Gaussent for agreeing to supervise this internship. It enabled me to discover the professional environment of a researcher, the timetable made up of numerous meetings and get-togethers, and the warm atmosphere of the laboratory. I also had the opportunity to attend seminars and a thesis defense, directed by Stéphane Gaussent. I'd like to thank him once again for the time he gave me and all he helped me discover.
	
	The internship ran from June 2 to July 11, 2025. I worked on my own, interspersed with meetings with my internship tutor. 
	I began by learning about the subject of my internship by reading the research articles Stéphane Gaussent had recommended. My first task was to understand them by making their content my own: I tried to redo the proofs or answer questions that came naturally (how do you properly define this object? How many diagrams are there?...). About twice a week, we'd get together to discuss what I'd understood, and Stéphane Gaussent would point me in new directions, suggest other articles or ask me new questions. I also had the opportunity to look up a few articles on my own, and a seminar helped me to progress in my research.
	
	\bigskip
	
	This report will give an almost chronological account of my findings, trying to highlight the difficulties I encountered before arriving at the final product.
	
	\bigskip
	
	In the section \ref{sec:algebre}, I'll start by defining Temperley-Lieb's diagrammatic algebra, as a space generated by diagrams that can be multiplied. Then I'll study its original version, as an algebra presented by generators and relations. Finally, I'll pose the key questions in my reasoning for the rest of the course, in particular the equivalence of definitions. The section \ref{sec:base} will be devoted to the search for a basis for algebra, playing on both planes: diagrammatic algebra and presented the algebra. The aim of this section is to create links between definitions, introduce the notion of rewriting and provide original proofs of certain results. I introduce a generalization of Temperley-Lieb algebra, in the section \ref{sec:orientee}, where this time the links are oriented. This will allow me to raise the need to introduce the category structure, on which the rewriting will be simpler. Finally, in the section \ref{sec:category}, I'll propose a presented category definition, and study the Temperley-Lieb category, using rewriting methods.

	\newpage
	\section{The Temperley-Lieb algebra}
	\label{sec:algebre}
	
	According to \cite{doty2024originstemperleyliebalgebraearly}, Temperley-Lieb algebra first appeared in \cite{2900521e-eda6-3524-a923-e614d35ed042} for the study of statistical mechanics problems. Baxter then formalized it as an presented algebra in \cite{alma990000435320306161} and Kauffman finally expressed it diagrammatically in \cite{HKAUFFMAN1987395}
	
	This algebra has been linked to many notions in mathematics, physics and computer science, and in particular in 2 articles I've looked into. Abramsky establishes links with knot theory and quantum mechanics in \cite{abramsky2009temperleyliebalgebraknottheory}, while Bowman et al. focuse on Kazhdan-Lusztig theory in \cite{bowman2024orientedtemperleyliebalgebrascombinatorial}.
	
	\begin{quote}
		\textbf{Aside.} On June 30, I had the opportunity to attend Paul Philippe's thesis defense on \textit{Kazhdan-Lusztig theory and masures}, which gave me a glimpse of what a thesis might look like. As the aim of the defense was to present his results in 45 minutes, I understood very little of what his work consisted of. My understanding is that it's a short presentation to a jury who know the subject and have read the thesis, but also to an open audience. This is followed by a series of questions from the members of the jury, which lasts as long as the presentation; they ask very precise questions that reveal the doctoral student's mastery of the subject.
	\end{quote}
	
	My point of view will be an abstract study of Temperley-Lieb algebra, with the aim of finding a basis for the algebra by applying rewriting theory. I will use Ridout and Saint-Aubin's pedagogical paper \cite{ridout2014standardmodulesinductiontemperleylieb} to develop the basic theory of our algebra.
	
	\subsection{Diagrammatic algebra}
	
	Although diagrammatic algebra came chronologically after presented algebra, it presents a more intuitive view of structure.
	
	\bigskip
	Let $n \in \mathbb{N}^*$.
	
	\begin{definition}
		A \textbf{n-diagram} is defined as a $\mathbb{R} ^2$ compact rectangle with two horizontal sides marked at regular intervals by $n$ points numbered from left to right, and where the points are connected 2 by 2 by a continuous curve, called a \textbf{link}. We also add the constraint that the links must be 2 by 2 disjoint.
	\end{definition}
	
	\begin{vocabulary} Here's a list of useful vocabulary for later on.
		\begin{itemize}
			\item A link connecting two points on opposite sides is called a \textbf{transverse link}.
			\item A transverse link connecting two points of the same index is called \textbf{straight link}.
			\item A link connecting two points on the same side is called \textbf{jump}.
			\item A jump connecting points of index $i$ and $i+1$ on the same line is called a simple jump of index $i$, where $i \in [\![ 1,n-1]\!]$.
			\item A \textbf{cap} is a single jump on the bottom line.
			\item A \textbf{cup} is a single jump on the top line.
		\end{itemize}
	\end{vocabulary}
	
	\begin{example} In the case $n=4$, below is a diagram consisting of a straight jump, a cap and a cup. The cap is a simple jump of index $3$.
		\begin{figure}[ht]
			\centering
			\figQuatreDiagramme
			\caption{4-diagram}
			\label{fig:4diag}
		\end{figure}
	\end{example}
	
	\FloatBarrier
	
	We identify the $n$-diagrams to within isotopy, i.e. we consider the quotient set $D_n$ of the $n$-diagrams by the isotopy equivalence relation. I won't go into detail on the notion of isotopy, my aim being just to have an idea to formalize the fact that two diagrams will be considered identical in our eyes if we can continuously deform, and without changing the anchor points, one to arrive at the other.
	
	\begin{definition}[$n$-Diagrammatic algebra] Let $\delta \in \mathbb{C}^{\times}$ \\
		We call \textbf{$n$-Diagrammatic algebra} and we denote $\mathcal{D}_n(\delta)$ the $\mathbb{C}$-vector space generated by $D_n$, provided with a multiplication defined for diagrams (and whose neutral is the diagram possessing only straight links) then extended by bilinearity to make $D_n$ an algebra : \\
		Let $d_1$ and $d_2$ be two $n$-diagrams, $d_1d_2$ is defined as $d_2$ is placed on top of $d_1$ where the top points of $d_1$ are identified with the bottom points of $d_2$ (then compressed to obtain the initial rectangle).
		To keep a diagram, we identify the \textbf{bubbles} with $\delta$.(See fig \ref{fig:delta})
	\end{definition} 
	
	\begin{figure}[ht]
		\centering
		\figDelta
		\caption{The multiplication $d_1d_2$ on the left results in the diagram (here $d_1$) with a straight link, a cup and a cap, multiplied by the scalar $\delta$.}
		\label{fig:delta}
	\end{figure}
	
	I started by counting the diagrams, for $n = 1,2,3 \dots$, but I couldn't see a formula for the general case. The article \cite{ridout2014standardmodulesinductiontemperleylieb} bijects diagrams with what I'm going to call \textbf{bridges}, which consist of the same thing as diagrams but on a line: we connect 2 by 2, without crossing, $2n$ points arranged on a line. The bijection consists in taking the upper segment and rotating it to send it to the lower line without creating a crossing. (See figure \ref{fig:dtop}, which repeats the example in figure \ref{fig:4diag}). This enabled me to find a recurrence formula using a direct method, i.e. without bijecting $D_n$ to the set of increasing paths under the diagonal, as we'll see later.
	
	\begin{figure}[ht]
		\centering
		\figDiagrammeToPont
		\caption{From a 4-diagram to an 8-bridge}
		\label{fig:dtop}
	\end{figure}
	
	\FloatBarrier
	
	\begin{theorem} Let $u_n$ be the number of n-diagrams, \\
		\[ u_0 = 1 \text{ and } u_n = \sum_{k}^{n-1}u_ku_{n-1-k} \; .\]
		\label{theorem:rec}
	\end{theorem}
	\begin{proof}
		We can see $u_0=1$ as "there's only one diagram with $0$ points, it's the empty rectangle".
		
		I'll count the number of bridges, noting that a bridge can only contain an even number of points so that each vertex is connected. \\ 
		To build a bridge with $2n$ points, I first construct an arc starting from the first point. Then there are $n$ possibilities for the other end of the arc, considering that a bridge must contain an even number of points. \\ 
		For each possibility, I can build a $2k$-point bridge to the left of the end of my arc and a $2(n-1-k)$-point bridge to the right, for a reasonable $k$: $k \in [\![ 0,n-1]\!]$. There are therefore $u_ku_{n-1-k}$ bridges possible in this configuration.
		
		All that's left to do is sum to obtain the recurrence formula.
	\end{proof}
	
	\begin{figure}[ht]
		\centering
		\figPont
		\caption{To the left of the $k=6$ endpoint, I can build a $2 \times 2$-point bridge and $2 \times 1$-point bridge to the right.}
		\label{fig:pont}
	\end{figure}
	
	\begin{remark}
	At this stage, I don't have a general formula for $u_n$ (depending only on $n$), but this property fully determines (by recurrence) a sequence. So I've decided not to dwell on it any further, since it will suffice later to check that the number of increasing paths under the diagonal, called the “Catalan number”, verifies \ref{theorem:rec}.
	
	What's more, I only want a few justifications of the soundness of the theory I'm about to study, and would rather spend more time working on Temperley-Lieb Algebra than justifying its existence.
	\end{remark}
	
	\subsection{Presented algebra}
	We can think of algebra as finding a minimal number of diagrams (generators) from which all the others can be obtained by multiplying them. The aim is then to write a diagram as a word from the generators. The next step is to find the minimum number of relations on these words in order to maintain isomorphism with the diagrammatic algebra. For the formal definition, I'll rely on one of my readings about braid theory : \cite{alma990005876670306161}.
	
	I'd like to point out that, originally, the opposite process took place.
	
	As I pointed out earlier, my aim is not to find these generators and relations, but to convince myself that they are indeed sufficient.
	
	\begin{notation} Let $n \in \mathbb{N}^*$, $i \in [\![1,n-1]\!]$,
		\begin{itemize}
			\item We note $\mathbf{1}$ or $\mathbf{1_n}$ the $n$-diagram with straight links only.
			\item We note $\mathbf{E_i}$ the $n$-diagram with a cup and a cap at index $i$, straight links elsewhere.
		\end{itemize}
	\end{notation}
	
	\begin{figure}[ht]
		\centering
		\figGenerateurs
		\caption{Respectively $1_4$, $E_1$, $E_2$, $E_3$}
		\label{fig:gen}
	\end{figure}
	
	These diagrams show three interesting relationships : \\
	$\forall i,j \in [\![1,n-1]\!]$
	\begin{itemize}
		\item $E_i^2=\delta E_i$ (already seen in figure \ref{fig:delta})
		\item If $|i-j|>1$, $E_iE_j=E_jE_i$ (see figure \ref{fig:com})
		\item $E_iE_{i \pm 1}E_i=E_i$ (see figure \ref{fig:pro})
	\end{itemize}
	
	\begin{figure}[ht]
		\centering
		\figCommutativite
		\caption{$E_3E_1=E_1E_3$}
		\label{fig:com}
	\end{figure}
	
	\begin{figure}[ht]
		\centering
		\figProche
		\caption{$E_1E_2E_1=E_1$}
		\label{fig:pro}
	\end{figure}
	
	\FloatBarrier
	
	\begin{definition}
		Let $S$ be a set, and let $S^n=S\times \dots \times S$ denote the set of n-uplets of elements of S. We define 
		\[S^* = \bigcup_{n\in \mathbb{N}} S^n\]
		the set of finite sequences of elements of S.
		
		A finite sequence of elements of S is called \textbf{word}, and $s_1s_2\dots s_n := (s_1,s_2,\dots,s_n)$.
		The empty word is denoted $\varnothing$.
		
		We define an internal law $\cdot$ on $S^*$, called \textbf{concatenation}, by 
		\[s_1s_2\dots s_n \cdot s_1's_2'\dots s_m' := s_1s_2\dots s_n s_1's_2'\dots s_m'\]
		This makes the notation for words consistent. This law is clearly associative, and admits the empty word as neutral, thus defining a monoid.
		
		The monoid $(S^*,\cdot)$ is called \textbf{free monoid}.
	\end{definition}
	
	\begin{definition} Let $S$ be a finite set of symbols. \\
		Let $A$ be the algebra generated by $S$, whose product is concatenation, and whose addition and multiplication by a scalar are formal.
		
		Let $R$ be a finite set of symmetric binary relations on $A$. \\
		Let $\mathcal{R}$ be the reflexive, transitive, linear closure of the relation generated by $R$ (the set $\{ (u'+uav+v',u'+ubv+v') | u,v,u',v' \in A, (a,b) \in R \} \cup \{ (u,u)|u\in A \}$)
		
		We define the \textbf{algebra presented by generators $S$ and relations $R$} as the quotient:
		\[<S|R> :=A/\mathcal{R}\]
	\end{definition}
	\begin{remark}
		I'm defining a structure presented in section \ref{sec:category} a little more neatly, which is how you'd have to go about it to get a complete definition. I'm content here with a compact definition.
		
		The quotient preserves the algebraic structure thanks to the definition of $\mathcal{R}$, and is well-defined because it is an equivalence relation.
		
		The algebra generated by a monoid $M$ (a group without the inverse axiom) is the set of finite linear combinations of elements of the monoid, where we define the product in such a way as to make it inevitably bilinear:
		
		\[ (\sum_{u \in U \subset S^*}\lambda_u u )*(\sum_{v \in V \subset S^*}\mu_v v) = \sum_{(u,v) \in U \times V \subset S^*}\lambda_u \mu_v uv \]
	\end{remark}
	
	\begin{definition}[Temperley-Lieb algebra] Let $\delta \in \mathbb{C}^{\times}$ \\
		We denote $TL_n(\delta)$ and we call \textbf{Temperley-Lieb algebra} the unitary algebra presented by the generators $e_1,\dots,e_{n-1}$ and the relations : \\
		$\forall i,j \in [\![1,n-1]\!]$
		\begin{itemize}
			\item $e_i^2=\delta e_i$
			\item If $|i-j|>1$, $e_ie_j=e_je_i$
			\item $e_ie_{i \pm 1}e_i=e_i$
		\end{itemize}
		Unitary means that the algebra has a neutral element, denoted by $\varnothing$.
	\end{definition}
	
	\subsection{Equivalence between definitions}
	
	It's time to convince ourselves that the two definitions of algebra are equivalent, i.e. that $TL_n(\delta)$ and $\mathcal{D}_n(\delta)$ are isomorphic.
	
	\begin{theorem}
		The application $\varphi : TL_n(\delta) \rightarrow \mathcal{D}_n(\delta)$ which associates $e_i$ generators with $E_i$ diagrams defines an algebra isomorphism.
		\label{theorem:iso}
	\end{theorem}
	
	Proof plan:
	\begin{enumerate}
		\item $TL_n(\delta)$ and $\mathcal{D}_n(\delta)$ are of the same finite dimension
		\item $\varphi$ is surjective
	\end{enumerate}
	Section \ref{sec:base} will give a proof of this theorem.
	
	\begin{remark}[Note on point 1]
		Our two algebras can be seen as algebras generated by a monoid :
		\begin{itemize}
			\item $\mathcal{D}_n(\delta)$ is generated by the monoid $D_n \cup \{\mathcal{T}\}$ where $\mathcal{T}$ commutes with all diagrams, $\mathcal{T}$ replaces $\delta$ when a bubble appears during multiplication, and with the identification $\mathcal{T} = \delta \cdot 1$
			\item $TL_n(\delta)$ is generated by the presented monoid $M_n = <e_1,\dots,e_{n-1}, \tau | R>$ where $e_i^2= \delta \cdot e_i$ is replaced by $e_i^2= \tau e_i$ with moreover $e_i \tau = \tau e_i$, then the identification $\tau = \delta \cdot \varnothing$.
		\end{itemize}
		Since $\mathcal{T}$ and $\tau$ are collinear with $1$ and $\varnothing$, we need to show that these monoids have the same cardinal, forgetting about $\mathcal{T}$ and $\tau$. To do this, \cite{ridout2014standardmodulesinductiontemperleylieb} puts them in bijection with paths growing under the diagonal, forgetting the monoid structure.
	\end{remark}
	
	\begin{definition} Let $n \in \mathbb{N}^*$.
		A $\mathbb{Z}^2$ path running from $(0,0)$ to $(n,n)$ with a $1$ step only up or to the right, forbidding crossing the diagonal, is called \textbf{increasing path under the diagonal}. 
		The set of paths growing under the diagonal is denoted $CCSD(n)$.
	\end{definition} 
	
	\begin{proposition}
		$D_n$ is in bijection with $CCSD(n)$.
	\end{proposition}
	\begin{proof}
		In fact, we put bridges in bijections with $CCSD(n)$.
		
		We read a bridge from left to right, taking a step to the right if a link opens and a step-up if a link closes. This is a path under the diagonal, because at any given moment there are more links open than closed.
		
		The reciprocal bijection consists in opening every link with every step to the right and closing the last open link with every step upwards, which is possible because the path is below the diagonal.
	\end{proof}
	
	\begin{remark}
		The next step is to use Jones' normal form to establish a bijection between $CCSD(n)$ and the presented monoid $M_n/{\sim} \cong <e_1,\dots,e_{n-1}| e_i^2=e_i,e_ie_j=e_je_i,e_ie_{i \pm 1}e_i=e_i>$. Where $\sim$ is the equivalence relation where two words are equal “to within $\tau$” ($v\sim w \Leftrightarrow \exists p \in \mathbb{N},v= \tau^p w \text{ or } w= \tau^p v$).
		
		However, this bijection does not give a word corresponding to a diagram. This means that the following \ref{fig:dnc} diagram is not commutative, as shown by \ref{fig:CCSDEi} and \ref{fig:CCSDei}.
		
		This gives a lot of interest to the algorithm I found in \ref{sec:FNJdiag}, which gives a word associated with a diagram, and which is also in normal form. So I'm going to abandon the \cite{ridout2014standardmodulesinductiontemperleylieb} approach using the $CCSD(n)$ set, to prove the \ref{theorem:iso} theorem using only Jones normal form.
	\end{remark}
	
	\begin{figure}[ht]
		\centering
		\figDiagNonCommutatif
		\caption{Non-commutative diagram}
		\label{fig:dnc}
	\end{figure}
	
	\begin{figure}[ht]
		\centering
		\figCCSDEi
		\caption{Image of $E_1$ by the bijection $D_n \rightarrow ponts(2n) \rightarrow CCSD(n)$}
		\label{fig:CCSDEi}
	\end{figure}
	
	\begin{figure}[ht]
		\centering
		\figCCSDei
		\caption{Image of $e_1$ by the bijection $M_n/{\sim} \space \rightarrow \space CCSD(n)$}
		\label{fig:CCSDei}
	\end{figure}

	\FloatBarrier
	\section{Basis search}
	\label{sec:base}
	
	Now that we've defined the Temperley-Lieb Algebra, several questions arise. Ridout and Saint-Aubin wonder (in \cite{ridout2014standardmodulesinductiontemperleylieb}) whether there is a representation of the algebra and what it is, while Abramsky studies the associated category (in \cite{abramsky2009temperleyliebalgebraknottheory}), to which I'll turn in section \ref{sec:category}. For my part, I'm looking for a basis for this algebra, which will enable me to introduce the notion at the heart of my internship: rewriting theory, based in particular on Philippe Malbos's detailed course \cite{MalbosPhilippe2019LoAR}.
	
	\subsection{Jones Normal Form}
	
	Here, our two algebras (not yet proved isomorphic) are each generated by a monoid. So we can find a basis by showing that the monoid is finite and exhibiting its elements. To do this, I'll reduce a word or diagram to a Jones normal form, and the elements of the monoid will then be the normal forms.
	
	\subsubsection{JNF of a word}
	
	Remember that $M_n$ is the monoid that generates the algebra $TL_n(\delta)$, and that it's a quotient. Exhibiting its elements therefore amounts to giving a system of representatives. I'm following Ridout and Saint-Aubin's method here.
	
	\begin{definition}
		A word $w \in S^*$ is in \textbf{Jones normal form} (JNF), if there exist $r \in \mathbb{N}^*$, $p \in \mathbb{N}$, $i_1<\dots<i_r$ and $j_1<\dots<j_r$ such that
		\[w = \tau^p(e_{i_1}e_{i_1-1} \dots e_{j_1}) \dots (e_{i_r}e_{i_r-1} \dots e_{j_r})\]
	\end{definition}
	
	\begin{definition}
		A word is said to be \textbf{reduced} if it is not possible to use relations to write it with fewer generators $e_1,\dots,e_{n-1}$ ($e_i^2$ is not reduced but $e_i\delta$ is).
		
		In other words, a word in $S^*$ is reduced if, as an element of $M_n$, it is not equal to a word written with fewer generators.
	\end{definition}
	
	\begin{lemma}
		The maximum index of a letter in a reduced word that is not a power of $\delta$ appears once.
	\end{lemma}
	\begin{proof}
		For $m \in \mathbb{N}^*$, we pose $H_m$ : \textit{"$\forall w \in M_n$, $w$ is reduced, is not a power of $\delta$, and is of maximum index m $\Rightarrow e_m$, appears once in w"}.
		
		\underline{For $m = 1$}. The only words with max index $1$ are products of $\tau$ and $e_1$, so $\tau^pe_i^q$ is written $q = 1$.
		
		\underline{Let $m>1$ be such that $H_1,\dots,H_{m-1}$ are true}. Let's show that $H_m$ is true. \\\
		Let $w \in M_n$ be reduced which is not a power of $\delta$, of maximum index $m$. Let's assume for the sake of argument that $e_m$ occurs at least twice.
		\[w = \dots e_mve_m \dots\]
		
		If $v$ is empty, we can simplify $e_m^2$, which is absurd.
		
		Otherwise, $m'<m$ is the maximum index of $v$. \\
		As $w$ is reduced, so is $v$, so $H_{m'}$ applies. Now, all the letters of $v$ except the single $e_{m'}$ are of index strictly less than $m-1$ and therefore commute with $e_m$. We can therefore write:
		\[w=\dots e_me_{m'}e_m \dots\]
		The two cases $m'=m-1$ and $m'<m-1$ result in a contradiction.
		Hence, $H_m$ is true.
		
		\underline{Conclusion.} By induction, $H_m$ is true for any non-zero integer $m$.
	\end{proof}
	
	\begin{theorem}
		The Jones normal forms form a system of representatives of $M_n$.
	\end{theorem}
	\begin{proof} We proceed as follows:
		\begin{enumerate}
			\item Every element of $M_n$ is equal to a reduced word.
			\item Any reduced word is equal as an element of $M_n$ to a word under FNJ.
		\end{enumerate}
		
		To avoid confusion between a word and an element of $M_n$ (equivalence class), I put a line above a word to indicate its equivalence class. \\
		No element is invertible in $M_n$ because no relation gives the empty word, so I'll omit the case of the empty word.
		
		\begin{par}{1.}
			Let $\overline{w}=\overline{e_{i_1} \dots e_{i_k}}$ be an element of $M_n$.
			
			The set of reduced words equal to $w$ in $M_n$ is included in the set of words on $\{e_1,\dots,e_{n-1}\}$ with less than $k$ letters, which is finite. So $w$ is equal to a reduced word $v$ in $M_n$ and is possibly itself.
		\end{par}
		
		\begin{par}{2.} 
			$v$ is reduced, so by the previous lemma its maximum index $m$ appears only once. We can then write $v = \dots e_m v'$ with $v'$ reduced.
			Switching $e_m$ to the right gives $\overline{v}=\overline{\dots e_me_{m-1} v''}$ and iterating this procedure $\overline{v}=\overline{u (e_me_{m-1} \dots e_{k})}$ with $k \leq m$.
			
			All that remains is to iterate the method over $u$, which is also reduced because $v$ is.
		\end{par}
	\end{proof}
	
	\begin{corollary}
		Since $\tau$ is collinear with $\varnothing$ the family of Jones normal forms, with no $\tau$ factor, generates $TL_n(\delta)$.
		\label{cor:familleGen}
	\end{corollary}
	
	\begin{remark}
		The proof of the theorem does not give an algorithm for obtaining the Jones normal form of a word.
	\end{remark}
	
	\subsubsection{JNF of a diagram}
	\label{sec:FNJdiag}
	Since we're looking for isomorphism between algebras, we want a diagram to be expressed from $E_i$, and we want the same properties to be verified. Thus, a diagram is supposed to be expressible in Jones normal form. So I decided to deviate from the method I'd been following up to now and look for an algorithm that gives the Jones normal form of a diagram. This consisted in adding capital letters to the FNJ for words, and replacing $\tau$ by $\mathcal{T}$.
	
	\begin{proposition}
		The Jones normal form is unique for diagrams, and so the family of factor-free $\mathcal{T}$ Jones normal forms is a free family of $\mathcal{D}_n(\delta)$.
		\label{prop:free}
	\end{proposition}
	\begin{proof} To convince yourself of this proof, I recommend visualizing $E_mE_{m-1} \dots E_{k}$ (see figure \ref{fig:segFNJ}).
		
		Suppose $d=mathcal{T}^p(E_{i_1}\dots E_{j_1})\dots (E_{i_r}\dots E_{j_r})=mathcal{T}^{p'}(E_{i_1'}\dots E_{j_1'})\dots (E_{i_{r'}‘}\dots E_{j_{r’}'})=d'$ are two FNJs.
		
		Since $j_{r-1} <j_r$, $E_{i_r}$ cannot simplify (there will always be $E_{i_r-2}$ between $E_{i_r-1}$,$E_{i_r}$ and $E_{i_r-1}$ ). This means that the link between the points at the top and bottom of index $i_r+1$ is not a straight link. Moreover, since $i_r$ is the maximum index of $d$, there are straight links at indices $i_r+2$ and above. Using the same arguments for $d'$, we obtain that $i_r = i_{r'}'$.
		
		Since $E_{j_r}$ is the last diagram of the product, it imposes a cup of index $j_r$ and none of higher index because again $j_{r-1} <j_r$. We deduce $j_r = j_{r'}‘$ using the same arguments on $d’$.
		
		By strict growth of the indices $i,j,i',j'$, no bubble can be created and therefore $p=p'$.
		
		Thus, two equal Jones normal forms are identical, by finite recurrence on the maximum index.
	\end{proof}
	
	\begin{figure}[ht]
		\centering
		\scalebox{1.2}{\figSegmentFNJ}
		\caption{$E_3E_2E_1$}
		\label{fig:segFNJ}
	\end{figure}
	
	\begin{algorithm}[ht]
		\caption{Jones normal form of a diagram.}\label{alg:FNJdiag}
		\begin{algorithmic}[1]
			\Require A diagram $d$
			\Ensure Jones normal form of the diagram $d$.
			\Function{JNF}{$d$}
			\If{$d=1_n$}
			\State \Return $1_n$
			\Else
			\State $j \gets$ index of the rightmost single jump $d$
			\State $i \gets -1+$ index of the last point without a right link
			\State $d$ is of the form :
			\State \figSeparationSegmentFNJ
			\State We write $d = d' E_iE_{i-1} \dots E_j$
			\State \Comment{\textit{Here, $d$ is arbitrary, and we've only identified the single jump of index $j$, shifted the links starting from positions $j+2$ to $i-3$ at the top and deformed the link starting from $i+1$ at the bottom.}}
			
			\State \Return \Call{JNF}{$d'$} $\cdot \ (E_iE_{i-1} \dots E_j)$
			\EndIf
			\EndFunction
		\end{algorithmic}
		\label{alg:FNJ}
	\end{algorithm}
	\FloatBarrier
	
	\begin{proposition}
		The previous algorithm is totally correct.
	\end{proposition}
	\begin{proof} First of all, the algorithm is well-defined. Indeed, $j$ exists because the only diagram with only cross-links is $1_n$, which is the case avoided by the first condition, and $j \leq i$ by definition. Moreover, $d'$ is well constructed, according to the commentary.
		\begin{itemize}
			\item[\underline{Ending}] $i$ is reduced by $1$ because if it is $0$, then $d=1_n$, is an integer and decreases strictly because the construction of $d'$ adds a right link at position $i+1$ and therefore the $i'$ calculated for the $d'$ diagram as input to the algorithm, will be less than i. Thus, the algorithm ends.
			\item[\underline{Correction}] The result is of the form $(E_{i_1}E_{i_1-1} \dots E_{j_1}) \dots (E_{i_r}E_{i_r-1} \dots E_{j_r})$ and we've seen that $i_1 < \dots < i_r$ ($i$ decreases but the concatenation is from the left).
			
			I now need to show that $j$ decreases strictly, i.e. that the position of the rightmost single jump in $d'$ is strictly less than $j$. There are two ways of creating a simple jump in $d'$ :
			\begin{itemize}
				\item If there were a jump above the simple jump of index $j$:
				
				\begin{tikzpicture}
					\node (1) at (0,0) {$j-1$};
					\node (2) [right=of 1] {$j$};
					\node (3) [right=of 2] {$j+1$};
					\node (4) [right=of 3] {$j+2$};
					\path[{Circle[]}-{Circle[]}] 
					(1) edge [bend right=90,looseness=1.7] (4)
					(2) edge [bend right=90,looseness=2] (3);
					\draw[dotted,thick] ($(1)+(0,-2)$) -- ($(4)+(0,-2)$);
				\end{tikzpicture}
				
				This simple jump has index $j-1<j$ in $d'$.
				\item If the deformation of the link from $i+1$ downwards creates a simple loop. In this case, the link from $i+1$ is transverse, and as it's not a straight link, by definition of $i$, and there must be no single jump after $j$, the link from $i+1$ down arrives before $j$ up.
			\end{itemize}
			Thus, $j_1 < \dots < j_r$, so the result is indeed in Jones normal form.
		\end{itemize}
		
		The algorithm is therefore totally correct.
	\end{proof}
	
	Since $\mathcal{T}$ is collinear with $1$, we deduce the following theorem from the proposition \ref{prop:free} and the algorithm \ref{alg:FNJ}.:
	\begin{theorem}
		The family of Jones normal forms, without factor $\mathcal{T}$, is a basis of $\mathcal{D}_n(\delta)$.
		\label{theorem:baseDiag}
	\end{theorem}
	
	\begin{par}{\begin{center}\textbf{Isomorphism between presented algebra and diagrammatic algebra}\end{center}}
		\begin{itemize}
			\item The image of literal Jones normal forms, without factor $\tau$, by the morphism $\varphi$ defined in theorem \ref{theorem:iso}, is a basis of the diagrammatic algebra according to theorem \ref{theorem:baseDiag}.
			\item The family of literal Jones normal forms, without factor $\tau$, is a generating family of the presented algebra according to corollary \ref{cor:familleGen}.
		\end{itemize}
		So $\varphi : TL_n(\delta) \rightarrow \mathcal{D}_n(\delta)$ is an algebra isomorphism.
	\end{par}
	
	\subsection{Rewriting theory}
	
	The aim of this subsection is to introduce an abstract notion: the rewriting system. Once proven convergent, it proves the existence of normal forms and automatically gives an algorithm for obtaining them. The only task will be to compare these normal forms with Jones'.
	
	\subsubsection{Rewriting systems}
	
	To introduce the notion of a rewriting system, I rely on Philippe Malbos's course \cite{MalbosPhilippe2019LoAR} and refer to it for proofs of results.
	
	\begin{definition}
		A set $\mathcal{A}$ and a family $(\rightarrow_i)_{i\in I}$ of binary relations on $\mathcal{A}$, called \textbf{rewriting relations}, are called \textbf{abstract rewriting system}.
	\end{definition}
	
	\begin{vocabulary} Let $a,b\in \mathcal{A}$
		\begin{itemize}
			\item $a \to b$ is called \textbf{rewriting step}.
			\item A \textbf{rewriting sequence} is a finite or infinite sequence of rewriting steps $a_1 \to a_2 \to a_3 \to \dots$.
			\item In the case of a finite sequence $a=a_1 \to a_2 \to \dots \to a_n =b$, we say that $a$ \textbf{is rewritten as} $b$ and we note $a\twoheadrightarrow b$ (the sequence is possibly of length 0: $a\twoheadrightarrow a$).
		\end{itemize}    
	\end{vocabulary}
	
	\begin{definition} Let $(\mathcal{A},\to_I)$ be a rewriting system.
		\begin{itemize}
			\item The system is said to be \textbf{terminating} if for any element $a \in \mathcal{A}$, we have that any sequence of rewrites starting with $a$ is finite.
			\item The system is said to be \textbf{confluent} if for each starting point $a \in \mathcal{A}$ such that $a \twoheadrightarrow b$ and $a \twoheadrightarrow c$, there exists $d \in \mathcal{A}$ such that $b \twoheadrightarrow d$ and $c \twoheadrightarrow d$ (see figure \ref{fig:conf}).
			\item The system is said to be \textbf{locally confluent} if for each starting point $a \in \mathcal{A}$ such that $a \to b$ and $a \to c$, there exists $d \in \mathcal{A}$ such that $b \twoheadrightarrow d$ and $c \twoheadrightarrow d$ (see figure \ref{fig:confloc}).
			\item The system is said to be \textbf{convergent} if it is terminating and confluent.
		\end{itemize}
	\end{definition}
	
	\begin{figure}[ht]
		\begin{subfigure}[c]{0.45\textwidth}
			\centering
			\figConfluence
			\caption{Confluence diagram}
			\label{fig:conf}
		\end{subfigure}
		\hfill
		\begin{subfigure}[c]{0.45\textwidth}
			\centering
			\figConfluenceLocale
			\caption{Local confluence diagram}
			\label{fig:confloc}
		\end{subfigure}
	\end{figure}
	
	\begin{definition} Let $a \in \mathcal{A}$,
		$a$ is said to be in \textbf{normal form} if no rewriting step starts from $a$.
	\end{definition}
	
	\begin{proposition}
		\begin{itemize}
			\item If a rewriting system is \underline{confluent} then every element admits at most one normal form.
			\item If a rewriting system is \underline{terminating} then every element admits at least one normal form.
		\end{itemize}
		Thus, if a rewriting system is \underline{convergent} then every element admits a unique normal form.
	\end{proposition}
	
	\begin{theorem}[Newman's lemma]
		If a rewriting system is \underline{locally confluent} and \underline{terminating} then it is convergent.
		\label{lemma:Newman}
	\end{theorem}
	
	We now consider rewriting systems whose elements are words in a given alphabet, and whose rewriting rules verify that if a word $a$ rewrites as a word $b$, then for all words $u$ and $v$, $uav$ rewrites as $ubv$.
	
	\begin{definition}
		We can compare local branches with the order \textbf{$\sqsubseteq$} generated by the relations
		\[(a,b,c)\sqsubseteq (uav,ubv,ucv)\]
		for any local branch (a,b,c) and all words u and v.
		
		An \textbf{overlap} is a branch that is neither of the form $(a,b,b)$ nor of the form $(ab,a'b,ab')$.
		
		A critical pair is a minimal overlap for the order $\sqsubseteq$.
	\end{definition}
	
	\begin{theorem}
		If all critical pairs of a rewriting system are confluent, then the system is locally confluent.
	\end{theorem}
	
	\subsubsection{Convergence for TL algebra}
	
	I'm going to look for a convergent rewriting system for the (presented) Temperley-Lieb algebra, which will give me uniqueness of normal forms and an algorithm for obtaining them: just look for a rewriting rule to apply in the word, apply it and repeat the process. This is totally correct according to the preceding theorems.
	
	In our framework, we study the rewriting system that guides the relations of the presented algebra. These rewriting rules make sense because we are transforming two words that are equal in the monoid.
	
	\begin{conjecture}
		The following rewriting system is convergent: (with the same indices as in the algebra presentation)
		
		The set is that of words on $\{delta,e_1,\dots,e_{n-1}\}$ and the relations are :
		\begin{align}
			e_i \delta &\to \delta e_i\\
			e_i^2 &\to \delta e_i \\
			e_i e_{i\pm 1} e_i &\to e_i \tag{$3\pm$}\\
			e_i e_j &\to e_j e_i \text{ if } j<i-1 \tag{4}
		\end{align}
	\end{conjecture}
	
	\begin{lemma}
		The previous rewrite system is terminating.
		\label{lemma:term}
	\end{lemma}
	\begin{proof}
		We order $\delta < e_1 < \dots < e_{n-1}$ and provide our set of words with lexicographic order.
		
		In our system, every relation $a \to b$ verifies $b<_{lex} a$.
		However, there is no strictly decreasing sequence for lexicographic order, so the system is terminating.
	\end{proof}
	
	To show the convergence of our system, all we need to do is demonstrate local confluence according to Newman's lemma \ref{lemma:Newman}, for which we study critical pairs, i.e. the moments when two rewriting rules overlap. For example, it's not necessary to deal with the case $e_3e_1 e_5e_6e_5$ where rules (3) and (4) are applied, because “they switch” in this case.
	
	I'm going to present the critical pairs in diagram form, trying not to deal with similar cases, specifically :
	
	\begin{itemize}
		\item I'm dealing with the case where rule (1) overlaps with rule (2), the other cases for rule (1) are similar, either shift delta to the left first, or apply the other rule first.
		\item Rule (1) does not overlap with itself.
		\item Rule (2) overlaps similarly to rule (1) with the other rules and with itself, because the overlap is on $e_i$, which is the letter that appears after rewriting.
		\item I deal with the case where rule (3+) overlaps with rule (3-) on one letter, the other $\pm$ cases being similar, and I deal with the case where the overlap is on two letters.
		\item I deal with the case where rule (4) overlaps with itself 
		\item I deal with the case where (3) overlaps with (4).
	\end{itemize}
	
	\begin{figure}[ht]
		\begin{subfigure}[b]{0.45\textwidth}
			\centering
			\begin{tikzpicture}[node distance=1.5cm]
				\matrix (m) [matrix of nodes, row sep=0.5cm, column sep=1cm] 
				{
					& $e_ie_i\delta$& 					 \\
					$e_i\delta e_i$ &   			&					 \\
					&   			& $\delta e_i\delta$ \\
					$\delta e_i^2$	&   			& 					 \\
					& $\delta^2e_i$& 				  	 \\
				};
				
				\path[->]
				(m-1-2) edge node[above, sloped] {(1)} (m-2-1) 
				(m-2-1) edge node[below, sloped] {(1)} (m-4-1) 
				(m-4-1) edge node[below, sloped] {(2)} (m-5-2) 
				
				(m-1-2) edge node[above, sloped] {(2)} (m-3-3)  
				(m-3-3) edge node[below, sloped] {(1)} (m-5-2); 
			\end{tikzpicture}
			\caption{Overlapping rules (1) and (2)}
		\end{subfigure}
		\hfill
		\begin{subfigure}[b]{0.45\textwidth}
			\centering
			\begin{tikzpicture}[node distance=1.5cm]
				\matrix (m) [matrix of nodes, row sep=1.5cm, column sep=0.5cm] 
				{
					& $e_ie_{i+1}e_ie_{i-1}e_i$	&	 		 	\\
					$e_ie_{i-1}e_i$ 	&   						&$e_ie_{i+1}e_i$\\
					& $e_i$						& 				\\
				};
				
				\path[->]
				(m-1-2) edge node[above, sloped] {($3+$)} (m-2-1)
				(m-2-1) edge node[below, sloped] {($3-$)} (m-3-2);
				\path[->]
				(m-1-2) edge node[above, sloped] {($3-$)} (m-2-3)
				(m-2-3) edge node[below, sloped] {($3+$)} (m-3-2);
			\end{tikzpicture}
			\caption{One-letter overlap of rules (3+) and ($3-$)}
		\end{subfigure}
		\hfill
		\begin{subfigure}[b]{0.45\textwidth}
			\centering
			\begin{tikzpicture}[node distance=1.5cm]
				\matrix (m) [matrix of nodes, row sep=2cm, column sep=0.5cm] 
				{
					& $e_ie_{i+1}e_ie_{i+1}$& 		 	 \\
					$e_ie_{i+1}$ 		&   					&$e_ie_{i+1}$\\
				};
				
				\draw[->] (m-1-2) -- node[above, sloped] {($3+$)} (m-2-1);
				\draw[dashed] (m-2-1) -- node[above, sloped] {$=$} (m-2-3);
				\draw[->] (m-1-2) -- node[above, sloped] {($3-$)} (m-2-3);
			\end{tikzpicture}
			\caption{Two-letters overlap of rules (3+) and ($3-$)}
		\end{subfigure}
		\hfill
		\begin{subfigure}[b]{0.45\textwidth}
			\centering
			\begin{tikzpicture}[node distance=1.5cm]
				\matrix (m) [matrix of nodes, row sep=1cm, column sep=0.5cm] 
				{
					& $e_ie_je_k$	& 		 	\\
					$e_je_ie_k$ &   			&$e_ie_ke_j$\\
					$e_je_ke_i$	&				&$e_ke_ie_j$\\
					& $e_ke_je_i$	&			\\
				};
				
				\path[->] 
				(m-1-2) edge node[above, sloped] {(4)} (m-2-1) 
				(m-2-1) edge node[below, sloped] {(4)} (m-3-1) 
				(m-3-1) edge node[below, sloped] {(4)} (m-4-2);
				\path[->] 
				(m-1-2) edge node[above, sloped] {(4)} (m-2-3)
				(m-2-3) edge node[above, sloped] {(4)} (m-3-3)
				(m-3-3) edge node[below, sloped] {(4)} (m-4-2);
			\end{tikzpicture}
			\caption{Overlapping rule (4) with itself, $k<j-1<i-2$}
		\end{subfigure}
		\caption{Confluence of critical pairs}
		\label{fig:pairesCritiques}
	\end{figure}
	
	\FloatBarrier
	
	I've shown local confluence for all critical pairs except the overlap of rules (3) and (4). This is in fact the only case of non-confluence.
	
	\begin{figure}[ht]
		\centering
		\begin{tikzpicture}
			\matrix (m) [matrix of nodes, row sep=1cm, column sep=0.5cm] 
			{
				& $e_ie_{i-1}e_{i}e_{i-2}$	& 		 				\\
				$e_ie_{i-2}$&   						&$e_ie_{i-1}e_{i-2}e_i$ \\
				$e_{i-2}e_i$&							&						\\
			};
			
			\path[->]
			(m-1-2) edge node[above,sloped] {($3-$)} (m-2-1)
			edge node[above,sloped] {(4)} (m-2-3)
			(m-2-1) edge node[below,sloped] {(4)} (m-3-1)
			(m-2-1.west) edge [bend right=45,dotted] node[left] {seule possibilité} (m-3-1.west)
			(m-2-3) edge [red,dashed] node[below,sloped] {(5)} (m-3-1);
			\path[thick] 
			($(m-3-1.south)+(-0.5,0)$) edge ($(m-3-1.south)+(0.5,0)$)
			($(m-3-1.south)+(-0.2,-0.1)$) edge ($(m-3-1.south)+(0.2,-0.1)$)
			($(m-2-3.south)+(-0.5,0)$) edge ($(m-2-3.south)+(0.5,0)$)
			($(m-2-3.south)+(-0.2,-0.1)$) edge ($(m-2-3.south)+(0.2,-0.1)$);
		\end{tikzpicture} 
		\caption{Overlapping rules (3) et (4)}
		\label{fig:nonConf}
	\end{figure}
	
	\FloatBarrier
	
	Two different, but still equal, normal forms appear in the monoid. To remedy this, we can apply Knuth-Bendix procedure, i.e., add a rewriting rule, symbolized by the red arrow, and study the new critical pairs, with the aim of making the system confluent. I then add the rules corresponding to the pathological cases $e_ie_{i-1}e_{i}e_{i-2}$ and $e_{i+2}e_ie_{i+1}e_i$ :
	\begin{align}
		e_ie_{i-1}e_{i-2}e_i &\to e_{i-2}e_i \tag{5} \\
		e_ie_{i+2}e_{i+1}e_i &\to e_ie_{i+2} \tag{6}
	\end{align}
	
	\FloatBarrier
	
	The overlaps of rule (5) with rules (1), (2), (3), (5) and (6) are similar to those of rule (3), and rule (6) is treated in the same way as rule (5). The study of the critical pair (5)-(4) $e_ie_{i-1}e_{i-2}e_{i}e_j$ has three cases: $j<i-3$ (similar to the critical pair (5)-(1)), $j=i-2$ and $j=i-3$.
	
	\begin{figure}[ht]
		\begin{subfigure}[b]{0.45\textwidth}
			\centering
			\begin{tikzpicture}[node distance=1.5cm]
				\matrix (m) [matrix of nodes, row sep=1cm, column sep=0.5cm] 
				{
					& $e_ie_{i-1}e_{i-2}e_{i}e_{i-2}$	& 		 		  	\\
					$e_{i-2}e_ie_{i-2}$ &  				&$e_ie_{i-1}e_{i-2}^2e_i$ 		\\
					$e_{i-2}^2e_i$		&				&$e_ie_{i-1}\delta e_{i-2}e_i$ 	\\
					$\delta e_{i-2}e_i$ &				&$\delta e_ie_{i-1}e_{i-2}e_i$	\\
				};
				
				\path[->]
				(m-1-2) edge node[above,sloped] {(5)} (m-2-1)
				edge node[above,sloped] {(4)} (m-2-3)
				(m-2-1) edge node[below,sloped] {(4)} (m-3-1)
				(m-3-1) edge node[below,sloped] {(2)} (m-4-1)
				(m-2-3) edge node[above,sloped] {(2)} (m-3-3)
				(m-3-3) edge[->>] node[above,sloped] {(1)} (m-4-3)
				(m-4-3) edge node[above,sloped] {(5)} (m-4-1);
			\end{tikzpicture}
			\caption{Non-pathological case: overlap of rules (5) and (4)}
		\end{subfigure}
		\hfill
		\begin{subfigure}[b]{0.45\textwidth}
			\centering
			\begin{tikzpicture}[node distance=1.5cm]
				\matrix (m) [matrix of nodes, row sep=1.5cm, column sep=0.5cm] 
				{
					& $e_ie_{i-1}e_{i-2}e_{i}e_{i-3}$	& 		 			  \\
					$e_{i-2}e_ie_{i-3}$&   					&$e_ie_{i-1}e_{i-2}e_{i-3}e_i$\\
					$e_{i-2}e_{i-3}e_i$&					&			 				  \\
				};
				
				\path[->]
				(m-1-2) edge node[above,sloped] {(5)} (m-2-1)
				edge node[above,sloped] {(4)} (m-2-3)
				(m-2-1) edge node[below,sloped] {(4)} (m-3-1)
				(m-2-3) edge [red,dashed] (m-3-1);
				\path[thick] 
				($(m-3-1.south)+(-0.5,0)$) edge ($(m-3-1.south)+(0.5,0)$)
				($(m-3-1.south)+(-0.2,-0.1)$) edge ($(m-3-1.south)+(0.2,-0.1)$)
				($(m-2-3.south)+(-0.5,0)$) edge ($(m-2-3.south)+(0.5,0)$)
				($(m-2-3.south)+(-0.2,-0.1)$) edge ($(m-2-3.south)+(0.2,-0.1)$);
			\end{tikzpicture} 
			\caption{Pathological case: overlapping rules (5) and (4)}
		\end{subfigure}
		\caption{Confluence of new critical pairs}
		\label{fig:nvPairesCritiques}
	\end{figure}
	
	By repeating the method, I add a finite number of rewriting rules. In fact, the pathological case occurs for $e_ie_{i-1}\dots e_{i-k}e_ie_{i-k-1}$, but to add the letter $e_{i-k-1}$, it must exist and therefore $i-k-1 \geq 1 \Leftrightarrow k \leq i-2 \leq n-3$. I thus obtain a confluent rewriting system.
	
	\begin{theorem}
		The following rewriting system is convergent:
		
		The set is that of the words on $\{\delta,e_1,\dots,e_{n-1}\}$ and the relations are :
		
		$\forall i,j \in [\![1,n-1]\!],\forall k \in[\![2,n-2]\!]$
		\begin{align}
			e_i \delta &\to \delta e_i                                      \tag{1} \\
			e_i^2 &\to \delta e_i                                           \tag{2} \\
			e_i e_{i\pm 1} e_i &\to e_i                                 \tag{$3\pm$}\\
			e_i e_j &\to e_j e_i \text{ si } j<i-1                          \tag{4} \\
			e_ie_{i-1}e_{i-2}\dots e_{i-k}e_i &\to e_{i-2}\dots e_{i-k}e_i  \tag{5} \\
			e_ie_{i+k}\dots e_{i+2}e_{i+1}e_i &\to e_ie_{i+k}\dots e_{i+2}  \tag{6}
		\end{align}
	\end{theorem}
	\begin{proof}
		The system is locally confluent according to the preceding study, and is terminating according to the same proof as for the \ref{lemma:term} lemma, since the added rules strictly decrease the lexicographic order.
		I deduce from Newman's lemma \ref{lemma:Newman} that the system is convergent.
	\end{proof}

	\subsubsection{Comparison between normal form via rewriting and Jones normal form}
	It's natural to wonder whether the Jones normal forms coincide with the normal forms by rewriting, given their very similar names, in addition to wondering whether we've obtained two different bases.
	
	What's more, if there is a coincidence, then we obtain the shape of the normal forms, without needing to run a computer to find them, which will be able to process only a precise $n$. To show that the normal forms are Jones' normal forms, I'll give an algorithm that transforms a word into its normal form, using only the rewriting rules. To do this, I visualize the words with a graph where the vertical axis represents the letters, the horizontal axis represents the position of the letter in the word.
	
	\begin{figure}[ht]
		\centering
		\figGrapheMot
		\caption{Exemple : visualisation du mot $e_1e_3\delta e_1$}
		\label{fig:exGrapheMot}
	\end{figure}
	
	I've listed the visualizations of the rewriting rules in Appendix \ref{sec:visu}. Before giving the algorithm, I'd like to point out the important points of visualizing a Jones normal form. Indeed, a Jones normal form is visually a sequence of antidiagonals, with strictly increasing slopes between the extremities (green slopes below, figure \ref{fig:exGrapheFNJ}).
	
	\begin{figure}[ht]
		\centering
		\figGrapheFNJ
		\caption{Example: visualizing a Jones normal form}
		\label{fig:exGrapheFNJ}
	\end{figure}
	
	\begin{algorithm}[ht]
		\caption{Jones' normal form of a word.}\label{alg:FNJmot}
		\begin{algorithmic}[1]
			\Require $w$ a word
			\Ensure Jones' normal form of the word $w$ using only rewrite rules.
			\Function{JNF}{$w$}
			\If{$w=\delta^*$ a possibly empty power of $\delta$}
			\State \Return $w$
			\Else
			\State Let $w=e_{i_1}\dots e_{i_r}$
			\State $m \gets \text{max}\{i_j \mid 1\leq j\leq r\}$
			\algstore{p1}
		\end{algorithmic}
	\end{algorithm}
	\begin{algorithm}[ht]
		\begin{algorithmic}[1]
			\algrestore{p1}
			\State Let $w = \green{v_0} \red{e_m} \green{v_1} \red{e_m} \dots \red{e_m} \green{v_s}$ where $s$ is the number of letters $e_m$ in $w$
			\State \figW 
			\State $v_1'\gets$ \Call{JNF}{$v_1$}
			\State \figVprime
			\State Apply rule (1) $e_m\delta \to \delta e_m$ to switch the $\delta$ to the left, then rule (4) $e_me_j\to e_je_m$ to switch the $v_1'$ antidiagonals at a distance (in height) greater than or equal to 2 from $e_m$.
			\State $v_0' \gets v_0$ concatenated with $\delta$ and the antidiagonals of $v_1'$ that begin with $e_j$, with $j<m-1$.
			\State $v_1'‘\gets e_{m-1}e_{m-2}\dots e_{m-k}$ the possible last antidiagonal of $v_1’$ starting with $e_{m-1}$.
			\If{$s \geq 2$}
			\State \figSdeuxVseconde
			\If{$v_i''=\varnothing$}
			\State Apply rule (2) $e_m^2 \to \delta e_m$
			\State \Return \Call{JNF}{$\green{v_0' \delta} \red{e_m} \green{v_2} \red{e_m} \dots \red{e_m} \green{v_k}$}
			\ElsIf{$v_1''=e_{m-1}$}
			\State Apply rule ($3-$) $e_me_{m-1}e_m \to e_m$
			\State \Return \Call{JNF}{$\green{v_0'} \red{e_m} \green{v_2} \red{e_m} \dots \red{e_m} \green{v_k}$}
			\algstore{p2}
		\end{algorithmic}
	\end{algorithm}
	\begin{algorithm}[ht]
		\begin{algorithmic}[1]
			\algrestore{p2}
			\Else
			\State Apply rule (5) $e_me_{m-1}e_{m-2}\dots e_{m-k}e_m \to e_{m-2}\dots e_{m-k}e_m$
			\State $v_0''\gets v_0'e_{m-2}\dots e_{m-k}$
			\State \Return \Call{JNF}{$\green{v_0''} \red{e_m} \green{v_2} \red{e_m} \dots \red{e_m} \green{v_k}$}
			\EndIf
			\Else{ $s=1$}
			\State $v_0''\gets$ \Call{JNF}{$v_0'$}
			\State \figSunVseconde
			\If{$v_0''=\delta^*$ a possibly empty power of $\delta$}
			\State \Return $\green{v_0''}\red{e_m}\green{v_1''}$
			\Else
			\State $j\gets$ index of the last letter of $v_0''$
			\If{$j<k$} (strictly positive black slope)
			\State \Return $\green{v_0''}\red{e_m}\green{v_1''}$
			\ElsIf{$j=m-1$}
			\State Apply rule ($3+$) $e_{m-1}e_me_{m-1} \to e_m$
			\State \Return \Call{JNF}{$\green{v_0''}e_{m-2}e_{m-3}\dots e_{m-k}$}
			\Else
			\State Apply rule (6) $e_je_{j+(m-j)}\dots e_{j+2}e_{j+1}e_j \to e_je_{j+(m-j)}\dots e_{j+2}$
			\State $u_0\gets (e_me_{m-1}\dots e_{j+2})$ 
			\State \textcolor{violet}{$u_1 \gets e_{j-1}\dots e_{m-k}$}
			\State Apply rule (4) to switch $u_0$ and \textcolor{violet}{$u_1$}.
			\State \figSunU
			\State \Return \Call{JNF}{$\green{v_0''} \textcolor{violet}{u_1}$} $u_0$
			\EndIf
			\EndIf
			\EndIf
			\EndIf
			\EndFunction
		\end{algorithmic}
	\end{algorithm}
	
	\FloatBarrier
	
	\begin{theorem}
		The previous algorithm is totally correct, so the normal forms of the rewriting system are the Jones normal forms.
	\end{theorem}
	\begin{proof} \textcolor{white}{.}
		\begin{itemize}
			\item[\underline{Ending}] $(m,\text{ number of letters } e_m)$ decreases strictly for the lexicographic order on $\mathbb{N}^2$ and ends for $(0,0)$, so the algorithm ends.
			\item[\underline{Correction}] \textcolor{white}{.}
			\begin{itemize}
				\item[\textit{Stop condition.}] If $w=\delta^*$ then $w$ is the word returned and is in FNJ.
				\item[\textit{Recursive calls}]
				Each recursive call returned (apart from the last one) is launched on a word obtained from a sequence of rewrites starting with $w$, or from a call to the algorithm on a subword, and is returned alone. We can therefore consider that if the algorithm is correct on these words, then the word returned is indeed in Jones normal form and obtained from rewriting steps.
				
				For the last return, FNJ($\green{v_0''} \textcolor{violet}{u_1}$) $u_0$, the aim is to see that the word returned is indeed in Jones normal form, given that it is obtained from a sequence of rewrites starting with $w$, or from calls to the algorithm on a sub-word. I first notice that when the word $\green{v_0''} \textcolor{violet}{u_1}$ is called by the algorithm, it is processed by the case $s=1$ (line 26). This step therefore consists in repairing the normal form by bringing any negative (black) slope downwards. What we need to show is that the last two antidiagonals are in Jones normal form, i.e. the last antidiagonal of FNJ($\green{v_0''} \textcolor{violet}{u_1}$) and $u_0$ have strictly positive (green) slopes between their ends. We'll look at each case of the FNJ($\green{v_0''} \textcolor{violet}{u_1}$) call. To avoid confusion, I underline the new values.
				
				If $\underline{v_0'‘}=\delta^*$ or $\underline{j}<\underline{k}$, FNJ($\green{v_0’'} \textcolor{violet}{u_1}$) $u_0$ is clearly in Jones normal form.
			
				If $\underline{j}=\underline{m}-1$, the last two antidiagonals are in Jones normal form, as only the first term of the penultimate antidiagonal has been removed.
				
				Otherwise, the penultimate antidiagonal is cut in half, the lower end is sent to the left and we check that the upper end concatenated with $u_0$ is indeed a Jones normal form. Clearly, the upper ends form a positive slope. For the lower ends, of index $\underline{j}+2$ and $j+2$, we have $\underline{j} < j$ because $\textcolor{violet}{u_1}$ is in Jones normal form, so $\underline{j}+2 < j+2$.
				
				The algorithm is therefore correct.
			\end{itemize}
		\end{itemize}
	\end{proof}
	
	\newpage
	\section{Oriented Temperley-Lieb algebra}
	\label{sec:orientee}
	
	Oriented Temperley-Lieb algebra is a generalization of Temperley-Lieb algebra, which consists in orienting diagram links (see figure \ref{fig:4diagOrient}). This algebra also admits a presentation (which I won't justify) and I refer to the article \cite{bowman2024orientedtemperleyliebalgebrascombinatorial} for its results and definitions. I'm looking for a rewriting system to find a basis for this new algebra.
	
	\begin{figure}[ht]
		\centering
		\figQuatreDiagrammeOriente
		\caption{4-oriented diagram}
		\label{fig:4diagOrient}
	\end{figure}

	\subsection{Presentation of algebra}
	
	In this subsection, I briefly review some of the results that will be useful to me.
	
	\begin{definition}
		The set $S_W$ is defined as the set of transpositions $s_i = (i,i+1)$ that generate $S_n$, the symmetric group.
		
		For $k \in [\![1,n-1]\!]$, we define the set ${W_k}$ of words on $S_W$ which, seen as permutations of $S_n$, are of minimal length. We construct this set starting from $s_k$, then concatenate by $s_i$ on the left for each $i$ and keep the word if we obtain a new permutation (in the sense that we see concatenation as a composition), and so on.
	\end{definition}
	\begin{remark}
		$W_k$ is a system of representatives of the quotient $S_n/(S_k\times S_{n-k+1})$. What will later be noted as $1_\lambda w 1_\mu$ for $\lambda,\mu \in W_k$ and $w$ a word on the letters $e_i$ is represented by the diagram where the southern orientation is that corresponding to $\lambda$, the northern orientation is that corresponding to $\mu$, and $w$ corresponds to the diagram without orientation.
	\end{remark}
	
	\begin{definition} Let $q$ be an indeterminate.    
		The \textbf{oriented Temperley-Lieb algebra}, $TLO_n(q)$, is defined as the $\mathbb{Z}[q,q^{-1}]$-unitary algebra generated by the elements 
		\[\{1_\lambda \mid \lambda \in W_k\} \cup \{ e_i \mid i \in [\![1,n-1]\!]\}\]
		and the following relations.  Idempotency relations,
		\begin{equation}
			\label{idempotentrel}
			\begin{split} {\sf 1}= \textstyle\sum_{\lambda\in W_k} {\sf 1}_\lambda, 
				\quad  
				{\sf 1}_\lambda e_i {\sf 1}_\lambda= 0   \text{ if }  	\lambda s_i \notin W_k,
				%
				\\
				{\sf 1}_\lambda  1_\mu =\delta_{\lambda,\mu}1_\lambda ,  \quad
				{\sf 1}_\lambda e_i {\sf 1}_\mu= 0   \text{ if }  	\mu \not \in \{\lambda,\lambda s_i \} 
			\end{split}
		\end{equation}
		for all $\lambda, \mu \in W_k$.
		for all $s_i\in S_W$, $\lambda\in W_k$, $\lambda s_i  \in W_k$ and $\mu, \nu \in \{\lambda, \lambda s_i\}$ :
		\begin{equation}\label{Esquare}
			{\sf 1}_\mu e_{ i}  {\sf 1}_\lambda   e_{ i}  {\sf 1}_\nu  =  
			q^{  \ell(\lambda s_i)-\ell(\lambda)} {\sf 1}_\mu e_i  {\sf 1}_\nu.
		\end{equation}
		If $|i-j|>1$, then
		\begin{equation}\label{commutation}
			e_i e_j=e_je_i, \qquad\quad 
			e_ie_{i\pm 1}e_i = e_i.
		\end{equation}
	\end{definition}
	
	\begin{remark}
		The relation ${\sf 1}= \textstyle\sum_{\lambda\in W_k}1_\lambda$ allows for more compact relations, but this is bad for rewriting. I build on Bowman et al.'s remarks in \cite{bowman2024orientedtemperleyliebalgebrascombinatorial}, that we can in fact replace generators with $1_\lambda$ and $1_\lambda e_i 1_\mu$, and insert $1_\lambda$ in the relations $1_\lambda e_i 1_\lambda e_j 1_\nu = 1_\lambda e_i 1_\nu e_j 1_\nu$ and $1_\mu e_i 1_\lambda e_{i+1} 1_\lambda e_i 1_\nu = 1_\mu e_i 1_\nu$. I refer to Bowman et al.'s article for quantifiers, as I won't need them, I'm just giving an idea here. Thus, an element of the algebra is a linear combination of words $w$ that can be written as follow :
		
		\[w = q^m \cdot (1_{\lambda_{(0,1)}} \dots 1_{\lambda_{(0,i_0)}}) e_1 (1_{\lambda_{(1,1)}} \dots 1_{\lambda_{(1,i_1)}}) \dots e_r (1_{\lambda_{(r,1)}} \dots 1_{\lambda_{(r,i_r)}})\]
		
		where the products of $1_{\lambda_{(i,j)}}$ simplify to 0 or $1_\lambda$. We therefore find words similar to those of non oriented algebra, but with 1s between the $e_i$:
		\begin{align*}
			e_1 \dots e_r &\in TL_n(\delta), \\
			1_{\lambda_0} e_1 1_{\lambda_1} \dots e_r 1_{\lambda_r}&\in TLO_{n,k}(q).
		\end{align*}
		
		\label{rem:rel}
	\end{remark}
	
	\subsection{Basis for normal forms}
	
	To provide a basis for oriented algebra, we can orient the relations of the remark \ref{rem:rel} and the other relations that have not been replaced. And see that the algebra is generated by a monoid, if we see $q$ as a letter. It is then a question of studying the rewriting system obtained by adding the necessary rules. I chose this rewriting system because it is very similar to the one described in the previous section, we would just have intercalated $1_\lambda$, for some $\lambda \in W_k$.
	\begin{center}
		\begin{math}
			\begin{array}{| r l | r l |}
				\multicolumn{4}{c}{\textbf{Comparaison between two word rewritting systems.}}\\
				\hline
				\multicolumn{2}{|c|}{\text{non oriented Temperley-Lieb}} & \multicolumn{2}{|c|}{\text{oriented Temperley-Lieb}} \\ 
				\hline
				\multicolumn{4}{c}{\text{Alphabet}}\\
				\hline
				\multicolumn{2}{|c|}{\{\delta,e_1,\dots,e_{n-1}\}}
				&
				\multicolumn{2}{|c|}{\{1_\lambda \mid \lambda \in {W_k}\} \cup \{q,1_{\lambda_1}e_11_{\mu_1},\dots,1_{\lambda_{n-1}}e_{n-1}1_{\mu_{n-1}}\}} \\
				\hline
				\multicolumn{4}{c}{\text{Rewritting rule}}\\
				\hline
				e_i \delta &\!\!\!\to\, \delta e_i 
				&
				e_i q \,\to\, q e_i \quad&\!\!\!\text{and}\quad 1_\lambda q \,\to\, q 1_\lambda
				\\			
				e_i^2 &\!\!\!\to\, \delta e_i
				&
				{\sf 1}_\mu e_{ i} {\sf 1}_\lambda e_{ i} {\sf 1}_\nu  &\!\!\!\to\, 
				q^{\ell(\lambda s_i)-\ell(\lambda)} {\sf 1}_\mu e_i {\sf 1}_\nu 
				\\
				e_i e_{i\pm 1} e_i &\!\!\!\to\, e_i
				&
				1_\mu e_i 1_\lambda e_{i\pm 1} 1_\lambda e_i 1_\nu &\!\!\!\to\, 1_\mu e_i 1_\nu
				\\			
				e_i e_j &\!\!\!\to\, e_j e_i \text{ if } j<i-1
				&
				1_\lambda e_i 1_\lambda e_j 1_\nu &\!\!\!\to\, 1_\lambda e_i 1_\nu e_j 1_\nu  \text{ if } j<i-1
				\\
				e_ie_{i-1}e_{i-2}\dots e_{i-k}e_i &\!\!\!\to\, e_{i-2}\dots e_{i-k}e_i
				&
				1_\mu e_i 1_{\lambda_0} e_{i-1} 1_{\lambda_1} \dots e_{i-k} 1_{\lambda_k} e_i 1_\nu &\!\!\!\to\, 1_{\lambda_1} e_{i-2}\dots e_{i-k} 1_{\lambda_k} e_i 1_\nu 
				\\
				e_ie_{i+k}\dots e_{i+2}e_{i+1}e_i &\!\!\!\to\, e_ie_{i+k}\dots e_{i+2}
				&
				1_\nu e_i 1_{\lambda_k} e_{i+k}\dots 1_{\lambda_1} e_{i+1} 1_{\lambda_0} e_i 1_\mu &\!\!\!\to\, 1_\nu e_i 1_{\lambda_k} e_{i+k}\dots e_{i+2} 1_{\lambda_1}
				\\ 
				&&
				{\sf 1}_\lambda e_i {\sf 1}_\lambda &\!\!\!\to\, 0 \text{ if }	\lambda s_i \notin {W_k}
				\\
				&&
				{\sf 1}_\lambda  1_\mu &\!\!\!\to\, \delta_{\lambda,\mu}1_\lambda
				\\
				&&
				{\sf 1}_\lambda e_i {\sf 1}_\mu &\!\!\!\to\, 0   \text{ if }  	\mu \not \in \{\lambda,\lambda s_i \} 
				\\
				\hline
			\end{array}
		\end{math}
	\end{center}
	
	For some $i,j,k,\lambda,\lambda_l,\mu,\mu_l$ and $\nu$.

	However, studying the confluence of the system seems tedious in notation but no more complicated. What's more, the article I'm relying on finds a basis (with another method than rewriting) that consists of Jones normal forms to which we add $1_\lambda$ to give orientation. It doesn't look like that, but on closer inspection, that's what it is. So I don't see the point of proving the result, with this method, or at least with this system of rewriting, because it's not going to do anything, except taint the theory of rewriting.
	
	That's why Stéphane Gaussent steered me towards the Temperley-Lieb category, which will make rewriting more visual and simple, after introducing complicated objects.

	\newpage
	\section{Temperley-Lieb category}
	\label{sec:category}
	
	In this section, I'm going to introduce the notion of category, which gives rise to a very important theory in all fields and can even be considered a foundation of mathematics. For the purposes of this course, I'll be confining myself to small categories (alive in set theory) and introducing only the notions I'll be using. My aim is to study (through the prism of rewriting) the Temperley-Lieb category whose set of endomorphisms of $n$ is isomorphic to the Temperley-Lieb algebra, in order to obtain a basis of Hom($n$,$m$) (sets of morphisms of $n$ in $m$).

	\subsection{Category theory}
	
	The general idea behind category theory is to move away from looking at interactions between elements of sets to looking at arrows between sets. For example, functions, group morphisms, linear applications... but above all forgetting about taking an element to return one.
	
	I'm going to rely on Mac Lane's book \cite{alma991007079110306161} for definitions. For my purposes, I see a category as an algebraic structure.
	
	\begin{definition}
		A \textbf{graph} is given by a set of objects $\cO$, arrows $\cF$ and two functions :
		\[\xymatrix{
			\cF \ar@<0.5ex>[r]^{\text{dom}} \ar@<-0.5ex>[r]_{\text{cod}} & \cO
		}\]
		
		Intuitively, dom$(f)$ (domain) represents the start of the arrow $f$ and cod$(f)$ (codomain) represents its arrival, so we note $f:\text{dom}(f)\to \text{cod}(f)$.
		
		We define the \textit{set of composable pairs} :
		\[\cF \times_{\cO} \cF = \{(g, f) \mid g,f \in \cF, \text{dom}(g) = \text{cod}(f) \}\]
	\end{definition}
	
	\begin{definition}
		A \textbf{category} $C$ is a graph with two additional functions
		\begin{equation*}
			\text{id} :
			\left\{
			\begin{array}{cll}
				\cO &\longrightarrow&\cF        \\
				a   &\longmapsto    &\text{id}_a
			\end{array}
			\right.
			,\quad
			\circ :
			\left\{
			\begin{array}{cll}
				\cF \times_{\cO} \cF&\longrightarrow& \cF \\
				(g, f)              &\longmapsto    & g \circ f
			\end{array}
			\right.
		\end{equation*}
		called \textit{identity} and \textit{composition}, such that
		\[
		\text{dom}(\text{id}_a) = a = \text{cod}(\text{id}_a)
		,\quad  \text{dom}(g \circ f) = \text{dom}(f)
		,\quad  \text{cod}(g \circ f) = \text{cod}(g)
		\]
		for all $a \in \cO$ and $(g,f) \in \cF \times_{\cO} \cF$. In addition, the following two axioms are verified:
		\begin{itemize}
			\item[]
			\item[] \textit{Associativity.} \quad $\forall (h,g,f) \in \cF \times_O \cF \times_O \cF, (h \circ g) \circ f = h \circ (g \circ f)$
			\item[]
			\item[] \textit{Identity law.}\quad $\forall f : a \to b, \text{id}_b \circ f = f = f \circ \text{id}_a$
		\end{itemize}
	\end{definition}
	
	\begin{definition}
		Arrows are sometimes called \textit{morphisms}, or \textit{homomorphism}.
		We can then define the set of homomorphisms$f : a \to b$
		\[\text{Hom}(a,b) = \{f \mid f\in\cF, \text{dom}(f) = a, \text{cod}(f) = b\}\]
		And the set of endomorphisms $f:a\to a$
		\[\text{End}(a) = \{f \mid f\in\cF, \text{dom}(f) = a = \text{cod}(f)\}\]
	\end{definition}
	
	\begin{definition}
		We call \textbf{strict monoidal category}, a category $C$ with an associative bifonctor $otimes:C \times C \to C$ :
		\[\otimes(\otimes \times 1) = (1\times \otimes)\otimes:C\times C\times C\to C\]
		and a neutral  $I$ :
		\[\otimes(I \times 1) = \text{id}_C = (1\times I)\otimes\]
	\end{definition}
	\begin{remark}
		I'll leave it to the reader to look up the definitions of bifonctor, $1$ and id$_C$ in \cite{alma991007079110306161}.
	\end{remark}
	
	\subsection{Strict presented linear monoidal category}
	
	I attended a seminar at the Camille Jordan Institute in Lyon on linear rewriting. The speaker was explaining how to obtain bases for the Hom$(a,b)$ vector spaces of a presented linear category, and as luck would have it, he introduced as an example the Temperley-Lieb category, which is in fact a simple case of a “diagrammatic category”. He noted:
	\[\cTL = <\figCap,\figCup \mid \figZigzag = \figId = \figZigzag[-1], \figBubble = 2 \varnothing>_{\otimes,\mathbb{C}}\]
	
	This helped me understand the appeal of this new object, on which rewriting can be done visually. He also talked about modulo rewriting, which consists in housing certain rules in the category structure and rewriting these rules “modulo”, i.e. without imposing meaning.
	
	While writing this report, I found a definition of category presented by generators and relations, more subtle than the one I propose below, in \cite{alma991007079110306161}. As I knew a definition of presented monoid, thanks to my reading of a book on braid theory \cite{alma990005876670306161}, I tried to adapt it for the notion of presented category. It's a rather crude vision of the structures presented by generators and relations, which requires a lot of work but which, I think, allows us to see the structure as a formal object, which only acts according to the rules we dictate to it. I hope to give an idea of the work that goes into an intuitive structure where an object is described by basic “bricks” and rules that make “parts of the wall fall apart”. However, I'm not giving the proofs here, which I haven't actually done, and it's possible that the way I build objects isn't complete, and still lacks precision. I'd like to have a better grasp of categories and Mac Lane's definition of a category, but I don't have the time during this course to look into it, so I'll leave this for my own future research.
	
	\begin{definition}
		Let $(M,\cdot_M)$ and $(N,\cdot_N)$ be two monoids. We define an internal law $\otimes$ on the Cartesian product $M\times N$ by
		\[(m,n)\otimes (m',n') := (m\cdot_M m',n\cdot_N n')\]
		This law is associative of neutral $(\varnothing_M,\varnothing_N)$ where $\varnothing_M$ is the neutral of $M$ and $\varnothing_N$ is the neutral of N.
		
		We then define the \textbf{product monoid} as the monoid $(M\times N,\otimes)$
	\end{definition}
	
	\begin{definition}
		Let $S,E$ be two sets and $F \subset E\times S\times S$ be a set of pairs of elements of $S$, labelled in $E$. Note $f:u\to v$ for $(f,u,v)\in F$.
		
		We define $F^\otimes$ the monoid generated by $F$ in the product monoid $E^*\times S^*\times S^*$.
		
		We define $(F^\otimes)^\circ$ the set of finite sequences $(f_n,u_n,v_n),\dots,(f_1,u_1,v_1)$ of elements of $F^\otimes$ such that $v_1=u_2, v_2=u_3,\dots,v_{n-1}=u_n$, and we associate a start $u_1$, an end $v_n$, and a label $f_n\circ\dots\circ f_1$.
		
		We then recursively define $((F^\otimes)^\circ)^{(n+1)}$ by $((F^\otimes)^\circ)^{(n)})^\otimes$ is the set of element sequences of $((F^\otimes)^\circ)^{(n)}$ whose start is the sequence of element starts, and similarly for the finish and label, and $((F^\otimes)^\circ)^{(n+1)} = ((((F^\otimes)^\circ)^{(n)})^\otimes)^\circ$ the set of sequences of elements of $(((F^\otimes)^\circ)^{(n)})^\otimes$ whose starts and finishes follow each other, and to which we associate the departure of the last element in the sequence, the arrival of the first element in the sequence, and the label of the sequence of labels of the elements in the sequence (noted with $\circ$).
		
		Finally, $F^{\otimes\circ*}=\bigcup_{n\in \mathbb{N}} ((F^\otimes)^\circ)^{(n)}$ is the set of alternating sequences of tensor products and compositions of arrows of F. We then denote $f:u\to v$ an element of this set, where $f$ is the label, $u$ is the start and $v$ is the finish.
		
		\bigskip
		
		Let $\cO=S^*$ and $\cF = F^{\otimes\circ*}$.
		
		We define the \textbf{free strict monoidal category}, denoted $<S\mid F>$ by the category whose objects are $\cO$ and arrows are $\cF$. The domain and codomain functions are defined by
		
		\begin{equation*}
			\forall f:u\to v \in \cF,\quad \text{dom}(f) = u,\quad \text{cod}(f) = v
		\end{equation*}
		
		The identity function is formally defined by
		\[\forall w \in \cO, id_w :w\to w\]
		and such that it verifies the identity law.
		Composition is defined as a finite sequence of two composable elements of $\cF$. This is clearly associative.
		
		We then define the bifonctor $\otimes$ as a sequence of 2 elements of $\cF$ for arrows and the concatenation of two words of $\cO$ for objects. The neutral is the empty word for objects and id$_\varnothing$ for arrows.
		
		Translated with DeepL.com (free version)
	\end{definition}
	\begin{remark}
		This defines a category, since the identity and composition functions verify the domain and codomain rules, and it should be shown that $\otimes$ is an associative bifonctor.
	\end{remark}
	\begin{notation}
		We also note $1_w$ for id$_w$.
		
		Since $\N$, the set of natural integers, is isomorphic to the free monoid on a 1-element set (we count the number of letters), we denote by integers the objects of a monoidal category presented by a generating object.
		
		Note $1$ for $1_1$, and $I$ for $id_\varnothing$.
	\end{notation}
	
	\begin{definition}
		Let $S,E$ be two sets, $F \subset E\times S\times S$ and $R\subset F^{\otimes\circ*}\times F^{\otimes\circ*}$ a set of symmetric binary relations on $F^{\otimes\circ*}$. The free strict monoidal category quotient by the equivalence relation generated by $R$, i.e. where the set of arrows is $\cF$ quotient by $R$, is called a \textbf{presented strict monoidal category}, which we denote $<S\mid F\mid R>$.
		
		A \textbf{strict monoidal linear category presented} consists in adding a linear structure.
	\end{definition}
	
	\begin{remark}
		It should be shown that the induced laws are well-defined.
		
		For an object $a\in\cO$, the set of endomorphisms of $a$ is an algebra.
	\end{remark}

	\subsection{Link to Temperley-Lieb algebra}
	
	We can now define the category of each of the two versions of the Temperley-Lieb algebra, oriented or non-oriented, in one fell swoop. All we have to do is change the generating objects.

	\begin{definition} Soit $\delta\in \mathbb{C}$\\
		The \textbf{Temperley-Lieb category} is defined as the strict monoidal linear category presented by a generating object and two generating arrows $cap : 2\to 0$ and $cup :0\to 2$ subject to the relations
		\begin{align*}
			(cap\otimes 1)\circ (1\otimes cup)=\ &1 =(1\otimes cap)\circ (cup\otimes 1)\\
			cap \circ cup &=\delta I
		\end{align*}
		with a special relationship, called the exchange relation,
		\[(1_p\circ f)\otimes(g\circ 1_m) = (f\circ 1_n)\otimes(1_q\circ g)\]
		for all arrows $f:n\to p$ and $g:m\to q$. We then abbreviate 
		\[\cTL(\delta) = <\figCap,\figCup \mid \figZigzag = \figId = \figZigzag[-1], \figBubble = \delta I,\textit{(exchange relation)}>\]
		We represent the exchange relation as follows:
		\[\figEchangeA = \figEchangeB\]
	\end{definition}
	
	\begin{proposition}
		End$_{\cTL(\delta)}(n)$ is isomorphic to $TL_n(\delta)$.
	\end{proposition}
	\begin{proof}
		The result is stated in \cite{abramsky2009temperleyliebalgebraknottheory}.
	\end{proof}
	
	In the same way, we define the Temperley-Lieb oriented category, and this time I'm only including the diagrams.
	
	\begin{definition} Let $q$ be an indeterminate.
		The \textbf{oriented Temperley-Lieb category}, denoted $\cTL\!\cO(q)$, is defined as the strict monoidal linear category presented by two generating objects
		\[\vee \text{ and } \wedge\]
		and four generating functions. $cap_+ : \vee\wedge \to \varnothing$, $cap_- : \wedge\vee \to \varnothing$, $cup_+ : \varnothing \to \vee\wedge$ and $cup_- : \varnothing \to \vee\wedge$
		\[ \figCapP \quad \figCapM \quad \figCupP \quad \figCupM \]
		subject to the relations 
		\begin{align*}
			\figZigzagUp = &\figIdUp = \figZigzagUp[-1]\\
			\figZigzagDown = &\figIdDown = \figZigzagDown[-1]\\
			\figBubbleTrigo &= q I \\
			\figBubbleHor &= q^{-1} I
		\end{align*}
		with exchange relation.
	\end{definition}
	
	\begin{remark}
		If we denote by $L_k$ the set of words on $\{\wedge,\vee\}$ containing exactly $k \; \vee$, everything suggests that the set of morphisms $d:w\to v \in \cTL\!\cO(q)$ where $w,v\in L_k$ is isomorphic to $TLO_{n,k}(q)$.
	\end{remark}

	\subsection{Rewriting for a basis of morphism algebras}
	
	The two previous definitions lead me to believe that for a good choice of equivalence relation and link between $\delta$ and $q$, it is possible to obtain that $\cTL(\delta)$ is a quotient $\cTL\!\cO(q)$. I'm therefore only studying a system for rewriting $\cTL\!\cO(q)$, which will enable me to obtain the bases of all sets of Hom$(v,w)$ morphisms (which are vector spaces), for $v,w \in \{\vee,\wedge\}^*$.
	
	\begin{theorem}
		The rewriting system on $\cTL\!\cO(q)$, with exchange relation, whose rules are :
		\begin{align}
			\figZigzagUp \to &\figIdUp \gets \figZigzagUp[-1] \tag{1-left/right}\\
			\figZigzagDown \to &\figIdDown \gets \figZigzagDown[-1] \tag{2-left/right}\\
			\figBubbleTrigo &\to q I 								\tag{3}\\
			\figBubbleHor &\to q^{-1} I								\tag{4}
		\end{align}
		is convergent.
	\end{theorem}
	\begin{remark}
		Modulo rewriting of certain relations means that we can use	these relations, in both directions of equality, before and after a rewriting step.
		
		We note $\to$ for a rewrite rule, $\rightsquigarrow$ or $\leftrightsquigarrow$ for a sequence  of relation which we work modulo and $\twoheadrightarrow$ for a possibly empty sequence of rewrite rules.
	\end{remark}
	
	\begin{figure}[ht]
		\centering
		\figConfluenceLocaleModulo
		\caption{Local confluence modulo}
		\label{fig:confMod}
	\end{figure}
	
	\begin{figure}[ht]
		\centering
		\figConfluenceLocaleModuloBis
		\caption{Equivalent version of local modulo confluence}
		\label{fig:confModBis}
	\end{figure}
	
	\FloatBarrier
	
	\begin{proof} \textcolor{white}{.}
		\begin{itemize}
			\item[\textit{Terminating.}]
			The number of generating functions (cap$_\pm$ and cup$_\pm$) decreases strictly for each rewrite rule and is unchanged by relations which we work modulo, so the rewrite system is convergent.
			\item[\textit{Confluent.}]
			There are two critical pairs, which differ only in their orientation, and are therefore treated in the same way. Below is a confluence diagram for one of the two critical pairs (figure \ref{fig:confCat}).
		\end{itemize}
		Newman's lemma for modulo rewriting (see \cite{MalbosPhilippe2019LoAR}) leads to the conclusion that the system converges.
	\end{proof}
	
	\begin{figure}[ht]
		\centering
		\figConfluenceCategorie
		\caption{Local modulo confluence diagram of a critical pair}
		\label{fig:confCat}
	\end{figure}
	
	Thus, for $v,w\in \{\vee,\wedge\}^*$, the set of normal forms $f:v\to w$ of the rewriting system forms a basis of the vector space of morphisms from $v$ into $w$, Hom$(v,w)$. Moreover, for the non-oriented case, the normal forms obtained are not the Jones normal forms of the diagrams, in the sense that the isomorphism which associates the endomorphism $1_{i-1} \otimes (cup\circ cap) \otimes 1_{n-i-1}$ with a generator $E_i$ of the diagram algebra, does not send a Jones normal form onto a normal form of the rewriting system of the category.
	
	\newpage
	\section*{Conclusion}
	
	The exploration of Temperley-Lieb algebra and rewriting theory involved refining methods to arrive at an efficient and explicit proof of the existence of a basis. I was able to make a small contribution, but above all I was able to discover some very interesting subjects, whether Temperley-Lieb algebra, rewriting theory or category theory. I was also able to rub shoulders with other fields from afar, by attending a seminar on fibrations, or Paul Philippe's thesis and discussions around his subject. This enabled me to get away from the usual classroom setting, without the pressure of having to master everything, and also to have time to devote to this.
	
	What's more, this internship has enabled me to strengthen my professional project. I'm attracted to mathematics, even if it's very theoretical, but I had little idea of what a researcher's job was, and a little more of what a teacher's job was. Stéphane Gaussent not only got me involved in a research activity, but also showed me many of the details that are part and parcel of a researcher's life: travel, meetings, small offices, supervising a thesis... I'm now all the more motivated to become a teacher-researcher, although I'm aware of the difficulties and drawbacks.
	
	Following this internship, I'll be able to take an interest in subjects I didn't even know existed until now. I intend to continue my studies at the magistère de mathématiques in Grenoble, while remaining curious.
	
	\newpage
	\printbibliography
	
	\newpage
	\appendix
	\section{Visualization of rewriting rules}
	\label{sec:visu}
	The vertical axis represents the letters, the horizontal axis the position in the word.
	
	\begin{figure}[ht]
		\centering
		\figRegleUn
		\caption{Rule (1) : $e_i\delta \to \delta e_i$}
		\label{fig:regleUn}
	\end{figure}
	
	\begin{figure}[ht]
		\centering
		\figRegleDeux
		\caption{Rule (2) : $e_i^2 \to \delta e_i$}
		\label{fig:regleDeux}
	\end{figure}
	
	\begin{figure}[ht]
		\centering
		\figRegleTroisPlus
		\caption{Rule (3+) : $e_ie_{i+1}e_i \to e_i$}
		\label{fig:regleTroisPlus}
	\end{figure}
	
	\begin{figure}[ht]
		\centering
		\figRegleTroisMoins
		\caption{Rule ($3-$) : $e_ie_{i-1}e_i \to e_i$}
		\label{fig:regleTroisMoins}
	\end{figure}
	
	\begin{figure}[ht]
		\centering
		\figRegleQuatre
		\caption{Rule (4) : $e_i e_j \to e_j e_i \text{ si } j<i-1$}
		\label{fig:regleQuatre}
	\end{figure}
	
	\begin{figure}[ht]
		\centering
		\figRegleCinq
		\caption{Rule (5) : $e_ie_{i-1}e_{i-2}\dots e_{i-k}e_i \to e_{i-2}\dots e_{i-k}e_i$}
		\label{fig:regleCinq}
	\end{figure}
	
	\begin{figure}[ht]
		\centering
		\figRegleSix
		\caption{Rule (6) : $e_ie_{i+k}\dots e_{i+2}e_{i+1}e_i \to e_ie_{i+k}\dots e_{i+2} $}
		\label{fig:regleSix}
	\end{figure}
\end{document}